\documentclass[]{interact}

\usepackage{epstopdf}
\usepackage[caption=false]{subfig}

\usepackage[numbers,sort&compress]{natbib}
\bibpunct[, ]{[}{]}{,}{n}{,}{,}
\renewcommand\bibfont{\fontsize{10}{12}\selectfont}

\theoremstyle{plain}
\newtheorem{theorem}{Theorem}[section]
\newtheorem{lemma}[theorem]{Lemma}
\newtheorem{corollary}[theorem]{Corollary}
\newtheorem{proposition}[theorem]{Proposition}

\theoremstyle{definition}

\theoremstyle{remark}
\newtheorem{remark}{Remark}[section]
\newtheorem{claim}{Claim}[section]

\usepackage{hyperref}
\usepackage{graphicx}
\usepackage{dsfont}
\usepackage{mathtools}
\usepackage{wrapfig}
\usepackage{amssymb, amsmath, latexsym}
\usepackage{nicefrac}
\usepackage{hhline}
\usepackage{multirow}
\usepackage{pifont}
\usepackage[colorinlistoftodos,bordercolor=orange,backgroundcolor=orange!20,linecolor=orange,textsize=scriptsize]{todonotes}
\usepackage{algorithm}\usepackage{algpseudocode}
\usepackage[symbol]{footmisc}

\usepackage{tikz}
\usetikzlibrary{matrix, positioning, fit}

\newcommand{\argmin}{\mathop{\arg\!\min}}

\def \RR {\mathbb R}

\newcommand{\EE}{\mathbf{E}}

\def\EE{\mathbb E}

\def\ag#1{{\color{black}#1}}

\def\pd#1{{\color{black}#1}} 

\begin{document}


%
\title{Accelerated gradient methods with absolute and relative noise in the gradient\thanks{
This work was supported by a grant for research centers in the field of artificial intelligence, provided by the Analytical Center for the Government of the Russian Federation in accordance with the subsidy agreement (agreement identifier 000000D730321P5Q0002) and the agreement with the Moscow Institute of Physics and Technology dated November 1, 2021 No. 70-2021-00138.
}}
%

%

\author{
\name{
Artem Vasin \textsuperscript{a} and Alexander Gasnikov \textsuperscript{a,b,c} and Pavel Dvurechensky\textsuperscript{d} and Vladimir Spokoiny\textsuperscript{d,e}
}
\affil{
\textsuperscript{a}Moscow Institute of Physics and Technology, Dolgoprudny, Russia;
\textsuperscript{b}Institute for Information Transmission Problems, Moscow, Russia;
\textsuperscript{c}ISP RAS Research Center for Trusted Artificial Intelligence, Russia;
\textsuperscript{d}Weierstrass Institute for Applied Analysis and Stochastics, Berlin, Germany;
\textsuperscript{e}Humboldt-University of Berlin, Berlin, Germany.
}
}

\date{Received: date / Accepted: date}
\maketitle              
\begin{abstract}
In this paper, we investigate accelerated first-order methods for smooth
convex optimization problems under inexact information on the gradient
of the objective. The noise
in the gradient is considered to be additive with two possibilities:
absolute noise bounded by a constant, and relative noise proportional
to the norm of the gradient. We investigate the accumulation of the errors
in the convex and strongly convex settings with the main difference with
most of the previous works being that the feasible set can be unbounded. \pd{The key to the latter is to prove a bound on the trajectory of the algorithm.} \pd{We also give a stopping criterion for the algorithm and consider extensions to the cases of stochastic optimization and composite nonsmooth problems.}
\end{abstract}


\section{Introduction}\label{Introduction}
We consider convex optimization problem on a closed convex  (not necessarily bounded) set \( Q \subseteq \mathbb{R}^n \):
\begin{equation}
\label{eq:pr_st}
    \min\limits_{x \in Q} f(x).
\end{equation}
We assume that the objective $f$ is $L_f$-smooth and strongly convex with the parameter $\mu \geqslant 0$, i.e., for all  $x,y\in Q$:
$$\|\nabla f(y) - \nabla f(x)\|_2 \leqslant L_f\|y - x\|_2,$$
$$f(x)+\langle\nabla f(x),y-x\rangle +\frac{\mu}{2}\|y - x\|_2^2 \leqslant f(y).$$

In the convergence rate analysis of different first-order methods these assumptions are typically used in the form of an upper and lower quadratic bounds \cite{devolder2013exactness,bubeck2015convex,ben-tal2015lectures,gasnikov2016efficient,beck2017first,gasnikov2017modern,nesterov2018lectures,lan2020first,tyurin2020adaptive,dvurechensky2020numerical,stonyakin2020adaptive,dvurechensky2021first,d2021acceleration} for the objective:
\begin{align}
\label{Lmu}
&f(x)+\langle\nabla f(x),y-x\rangle +\frac{\mu}{2}\|y - x\|_2^2 \leqslant f(y) \nonumber \\ 
&\hspace{3cm}\leqslant f(x)+\langle\nabla f(x),y-x\rangle + \frac{L_f}{2}\|y-x\|_2^2.
\end{align}
Note that the last relation is a consequence of the $L_f$-smoothness and, in general, is not equivalent, to it \cite{taylor2017smooth,gasnikov2017modern}.

In many applications, instead of an access to the exact gradient $\nabla f(x)$ an algorithm has access only  to its inexact approximation $\tilde{\nabla} f(x)$. Typical examples include gradient-free (or zeroth-order) methods which use a gradient estimator based on finite differences \cite{conn2009introduction,risteski2016algorithms,berahas2021theoretical}, and optimization problems in infinite-dimensional spaces related to inverse problems \cite{kabanikhin2011inverse,gasnikov2017convex}.
The two most popular definitions of gradient inexactness in practice are \cite{polyak1987introduction} as follows: for all $x \in Q$ it holds that
\begin{equation}\label{inexact}
    \|\tilde{\nabla} f(x) - \nabla f(x)\|_2 \leqslant \delta, \text{ (absolute error) or }
\end{equation}
\begin{equation}\label{relative_inexact}
    \|\tilde{\nabla} f(x) - \nabla f(x)\|_2 \leqslant \alpha\|\nabla f(x)\|_2, \quad \alpha\in[0,1) \text{ (relative error)}.
\end{equation}
Under assumption \eqref{inexact}, many results exist for non-accelerated and accelerated first-order methods, see, e.g., \cite{polyak1987introduction,d2008smooth,cohen2018acceleration,ajalloeian2020analysis}. These results are in a sense pessimistic in general with the explanation going back to the analysis in \cite{poljak1981iterative}. We can explain this by a very simple example. Consider the following problem
\begin{equation}\label{example}
 \min_{x\in\RR^n} \left\{f(x)=\frac{1}{2}\sum_{i=1}^n \lambda_i \cdot(x^i)^2 \right\},   
\end{equation}
where $0\leqslant \mu = \lambda_1 \leqslant \lambda_2 \leqslant ...\leqslant \lambda_n = L_f$, $L_f\ge2\mu$. Clearly, the solution of this problem is $x^* = 0$. Assume that the inexactness takes place only in the first component $x_1$, i.e., instead of $\partial f(x) / \partial x^1 = \mu x^1$ we have access to $\tilde{\partial} f(x) / \partial x^1 = \mu x^1 - \delta$, where $\delta$ is the error. For the simple gradient descent
$$x_{k} = x_{k-1} - \frac{1}{L_f}\tilde{\nabla}f(x_{k-1}),$$
we can conclude that if $x_0^1\geqslant 0$, then for all $k\in\mathbb{N}$ large enough, i.e., $k\gg L/\mu$, it holds that
\begin{equation}\label{lowerbound}
 x_k^1 \geqslant \frac{\delta}{L}\frac{1 - (1 -\mu/L_f)^k}{1-(1 -\mu/L_f)}\simeq \frac{\delta}{\mu}. 
\end{equation}
Hence,\footnote{This bound corresponds to the worst-case philosophy, i.e., choosing the worst example for the considered class of methods \cite{nemirovsky1983problem,nesterov2018lectures,bubeck2015convex,gasnikov2017modern}. We expect more interesting results by considering average-case complexity  \cite{scieur2020universal,pedregosa2020average}.}
$$f(x_k) - f(x^*) \gtrsim \frac{\delta^2}{2\mu}.$$
From this result, we see that it may be problematic to approximate $f(x^*)$ with any desired accuracy, especially in the ill-conditioned setting when the strong convexity constant $\mu$ is smaller than the desired accuracy $\varepsilon$.
For accelerated gradient methods the situation may be even worse since they are more sensitive to the gradient errors and such errors may even be accumulated by the algorithm \cite{devolder2014first,dvurechensky2016stochastic,gasnikov2017modern}.
This drawback may be overcome by proposing a certain stopping rule so that the algorithm does not try to minimize below some threshold given by the gradient error or by adding a strongly convex regularizer with coefficient $\mu$ of the same order as the desired accuracy $\varepsilon$, see \cite{poljak1981iterative,polyak1987introduction,nemirovski1986regularizing,gasnikov2017modern}. 
Roughly speaking,  for non-accelerated algorithms it was proved in \cite{poljak1981iterative,polyak1987introduction} that if $\delta$ is of the order $\varepsilon^2$, then it is possible to reach $\varepsilon$-accuracy in the objective residual function in almost the same number of iterations as in the exact case $\delta = 0$ by applying a computationally convenient stopping rule.


In this paper, we analyze an accelerated gradient method in both convex and strongly convex settings and estimate how the gradient error defined in \eqref{inexact} influences the convergence rate. An important part of our contribution is that our analysis is made without an assumption that the feasible set $Q$ is bounded. \pd{The main key for this development is a recurrent estimate for the distance between the current iterates and the optimal solution closest to the starting point.} In particular, our results imply that it is sufficient to assume that $\delta$ is of the order $\varepsilon$ in order to obtain objective residual of the order $\varepsilon$. We also present a stopping rule and prove that if it is satisfied at some iteration, the algorithm solves problem \eqref{eq:pr_st} with certain accuracy. Moreover, we prove that until this rule is fulfilled, the trajectory of the algorithm is bounded \pd{(which helps us to treat the setting of possibly unbounded set $Q$)} and that it is fulfilled for sure in a number of iterations which is optimal for the class of smooth convex optimization problems.


Under assumption \eqref{relative_inexact}, non-accelerated gradient method for strongly convex problems is shown in \cite{polyak1987introduction} to have linear convergence with condition number $O\left(\frac{1}{1-\alpha}\cdot \frac{L_f}{\mu} \right)$, i.e. $\frac{1}{1-\alpha}$ times worse than in the exact case. Yet, convergence to any small error is guaranteed unlike the case of inexactness \eqref{inexact}. 
This result holds also under the relaxed strong convexity assumption \cite{gasnikov2017modern} known as Polyak--Lojasiewicz or gradient domination condition. We are not aware of any such results for accelerated gradient methods.


In this paper, we analyze an accelerated gradient method under inexact gradients satisfying \eqref{relative_inexact} and answer the question of what is the maximum value of $\alpha$ such that the accelerated algorithm with inexact gradients converges with the same rate as the exact accelerated algorithm. For the case $\mu \ne 0$ our answer is that $\alpha$ should satisfy $\alpha=O\left(\frac{\mu}{L_f}\right)$. We hypothesise that this bound can be improved to $\alpha=O\left( \left(\frac{\mu}{L_f}\right)^{\pd{1/2}}\right)$ and, for the case $\mu=0$, the iteration-dependent value $\alpha_k$ should satisfy $\alpha_k=O\left( \left(\frac{1}{k}\right)^{3/2}\right)$, where $k$ is the iteration counter.
Numerical experiments demonstrate that, in general, for $\alpha$ larger than the mentioned above thresholds the convergence may slow down a lot up to divergence for the considered accelerated method.


Close results with the bound $\alpha=O\left( \left(\frac{\mu}{L_f}\right)^{5/4}\right)$ in the case $\mu \gg \varepsilon$ were recently obtained using another techniques in stochastic optimization with decision dependent distribution \cite{drusvyatskiy2022stochastic} and policy evaluation in reinforcement learning via reduction to stochastic variational inequality with Markovian noise \cite{kotsalis2022simple}. In \cite{kotsalis2022simple,drusvyatskiy2022stochastic}, the authors assumed that
\begin{equation}\label{relative_inexactDL}
    \|\tilde{\nabla} f(x) - \nabla f(x)\|_2 \leqslant B\|x - x^*\|_2.
\end{equation}
Since $x^*$ is a solution, when $Q=\mathbb{R}^n$, we have $\nabla f(x^*) = 0$. 
Therefore, $$\|\nabla f(x) \|_2 = \|\nabla f(x) - \nabla f(x^*)\|_2 \leqslant L_f\|x - x^*\|_2.$$
Thus, if \eqref{relative_inexact} holds, then \eqref{relative_inexactDL} also holds with $B = \alpha L_f$.


\section{Ideas behind the results}\label{Ideas} 

\subsection{Absolute noise}
Important results on gradient error accumulation for first-order methods were developed in a series of works by O.~Devolder, F.~Glineur and Yu. Nesterov 2011--2014 \cite{devolder2011stochastic,devolder2014first,devolder2013first,devolder2013exactness}. In these works, the authors were motivated by inequalities \eqref{Lmu}. Their idea was to relax \eqref{Lmu}, assuming inexactness in the gradient, introducing the inexact gradient $\tilde{\nabla} f(x)$, satisfying  for all $x,y\in Q$
\begin{align}
\label{inexact_model_1}
&f(x)+\langle\tilde{\nabla} f(x),y-x\rangle +\frac{\mu}{2}\|y - x\|_2^2- \delta \leqslant f(y) \nonumber \\ 
&\hspace{3cm}\leqslant f(x)+\langle\tilde{\nabla} f(x),y-x\rangle + \frac{L_f}{2}\|y-x\|_2^2 +\delta.
\end{align}
This assumption allows to develop a theory for error accumulation for first-order methods. 
In particular,  they obtained the following convergence rates for non-accelerated gradient methods:
\begin{equation}\label{convergenceNA} 
f(x_k) - f(x^*) =O\left(\min\left\{ \frac{L_f R^2}{k} +  \delta,  L_f R^2\exp\left(-\frac{\mu}{L_f}k\right) + \delta\right\}\right), 
\end{equation}
 and for accelerated methods:
\begin{equation}\label{convergenceA} 
f(x_k) - f(x^*) =O\left(\min\left\{ \frac{L_f R^2}{k^2} +  k\delta,  L_f R^2\exp\left(-\sqrt{\frac{\mu}{L_f }}\frac{k}{2}\right) +  \sqrt{\frac{L_f }{\mu}}\delta\right\}\right), 
\end{equation} 
where $R$ is such that $ \|x_{start} - x^*\|_2 \leqslant R$, i.e., an estimate for the distance between the starting point $x_{start}$ and a solution $x^*$. If \pd{$x^{*}$} is not unique, one may take \pd{$x^{*}$} to be the closest point to $x_{start}$. 
Both of these bounds are unimprovable \cite{devolder2014first,devolder2013first}. See also \cite{devolder2013exactness,dvurechensky2016stochastic,kamzolov2020universal} for \pd{``intermediate''} situations between accelerated and non-accelerated methods and extensions for stochastic optimization.

Following \cite{devolder2013first}, it is possible to make a reduction of the ``absolute noise'' inexactness in the sense of \eqref{inexact} to  the inexactness in the sense of \eqref{inexact_model_1} by setting
\begin{equation}\label{deltamu}
\delta = \delta_{\eqref{inexact_model_1}} = \frac{\delta_{\eqref{inexact}}^2}{2L_f} + \frac{\delta_{\eqref{inexact}}^2}{\mu} \simeq \frac{\delta_{\eqref{inexact}}^2}{\mu}
\end{equation}
and setting $L_{f,\eqref{inexact_model_1}}=2L_{f,\eqref{Lmu}}$, $\mu_{f,\eqref{inexact_model_1}}=\mu_{f,\eqref{Lmu}}/2$. The key observations here are that
    $$\langle\tilde{\nabla} f(x) - \nabla f(x),y-x\rangle \leqslant \frac{1}{2L_f}\|\tilde{\nabla} f(x) - \nabla f(x)\|_2^2 + \frac{L_f}{2}\|y-x\|_2^2,$$
$$\langle\tilde{\nabla} f(x) - \nabla f(x),y-x\rangle \geqslant \frac{1}{\mu}\|\tilde{\nabla} f(x) - \nabla f(x)\|_2^2 - \frac{\mu}{4}\|y-x\|_2^2.$$
From this reduction, we see that when $\mu > 0$, for non-accelerated methods.  the result \eqref{convergenceNA} is almost the same as in the example in \eqref{example}. We see also, that, if the error can be controlled, to guarantee that $f(x_k) - f(x^*) \leqslant \varepsilon$ for non-accelerated method when
\footnote{If $\mu \lesssim \varepsilon$, we can regularize the problem and guarantee that $\mu = \Omega( \varepsilon)$, see \cite{gasnikov2017modern}. Another advantage of strong convexity is the possibility to use the norm of inexact gradient for the stopping criteria, see \cite{gasnikov2017modern,poljak1981iterative}. Yet, regularization requires \cite{gasnikov2017modern} some prior knowledge about the distance to the solution. Since we typically don not have such information the procedure becomes more difficult via applying the restart technique, see \cite{gasnikov2016efficient,gasnikov2017modern}.} 
$\mu = \Omega( \varepsilon)$ we should set $\delta_{\eqref{inexact}}= O(\varepsilon)$, which is an expected result. Unfortunately, for accelerated methods, such reduction leads to the bound  $\delta_{\eqref{inexact}} =O( \varepsilon^{3/2})$, which is worse than our bound indicated in  Section~\ref{Introduction}. The key to our improvement is a more refined version of \eqref{inexact_model_1}.

In the works \cite{devolder2014first,dvinskikh2019decentralized,dvinskikh2020accelerated,stonyakin2021inexact,stonyakin2020adaptive}, the following refined version of \eqref{inexact_model_1} is used:
\begin{align}
\label{inexact_model_2}
&f(x)+\langle\tilde{\nabla} f(x),y-x\rangle +\frac{\mu}{2}\|y - x\|_2^2- \delta_1\|y-x\|_2 \leqslant f(y) \nonumber \\ 
&\hspace{3cm}\leqslant f(x)+\langle\tilde{\nabla} f(x),y-x\rangle + \frac{L_f}{2}\|y-x\|_2^2 +\delta_2.
\end{align}
These inequalities lead to the following counterparts of \eqref{convergenceNA}  and \eqref{convergenceA}, respectively, for non-accelerated gradient methods: 
\begin{align}\label{convergenceDNA} 
&f(x_k) - f(x^*) \nonumber \\ 
&\hspace{0cm} =O\left(\min\left\{ \frac{L_f R^2}{k} + \tilde{R}\delta_1 + \delta_2,  L_f R^2\exp\left(-\frac{\mu}{L_f }k\right) + \tilde{R}\delta_1 + \delta_2\right\}\right), 
\end{align}
and for accelerated gradient methods \cite{devolder2014first,dvinskikh2020accelerated}:  
\begin{align}\label{convergenceDA} 
&f(x_k) - f(x^*) \nonumber \\ 
&\hspace{0cm}
=O\left(\min\left\{ \frac{L_f R^2}{k^2} + \tilde{R}\delta_1 + k\delta_2,  L_f R^2\exp\left(-\sqrt{\frac{\mu}{L_f }}\frac{k}{2}\right) +  \tilde{R}\delta_1 + \sqrt{\frac{L_f }{\mu}}\delta_2\right\}\right), 
\end{align}
where $\tilde{R}$ is the maximum distance between the sequences of iterates generated by the algorithm and the solution $x^*$ closest to the starting point.
 
From \eqref{convergenceDNA}, \eqref{convergenceDA}, we see that if $\tilde{R}$ is bounded,\footnote{In many situations this is true. For example, when $Q$ is bounded or when $\mu\gg\varepsilon$.} then by setting
 \begin{equation*}
    \delta_1 = \delta_{\eqref{inexact}}, \delta_2 = \frac{\delta_{\eqref{inexact}}^2}{2L_f }, 
 \end{equation*}
  we obtain the desired result: it is possible to guarantee $f(x_k) - f(x^*) \leqslant \varepsilon$ with $\delta_{\eqref{inexact}}=O(\varepsilon)$.  
 
Previous works mainly rely on the assumption that  $\tilde{R}$ is bounded. As we may see from example~\eqref{example}, in general, when the strong convexity parameter $\mu$ is small compared to the desired accuracy $\varepsilon$, only a bound 
$$\tilde{R}\simeq R + \frac{\delta_{\eqref{inexact}}}{\mu}\gtrsim R + \frac{\delta_{\eqref{inexact}}}{\varepsilon}$$ 
is possible to obtain \cite{gasnikov2017modern}. This bound leads to very pessimistic estimates. Moreover, the growth of $\tilde{R}$ is observed in different numerical experiments and in theoretical estimates caused by error accumulation. In our work, we investigate this problem and, in particular, propose an alternative to regularization\footnote{By using regularization we can guarantee $\mu\sim\varepsilon$ and therefore with $\delta_{\eqref{inexact}}\sim\varepsilon$ we have the desired estimate $\tilde{R}\simeq R$.} approach that is based on ``early stopping''\footnote{This terminology is popular also in Machine Learning community, where ``early stopping'' is used also as an alternative to regularization to prevent overfitting   \cite{goodfellow2016deep}.} of the considered iterative procedure by developing proper stopping rule.

\subsection{Relative noise}
We now explain a way of reduction of the relative inexactness in the sense of \eqref{relative_inexact} to the inexactness in the sense of \eqref{inexact_model_1}, which allows us to apply \eqref{convergenceA} when $\mu\gg\varepsilon$. Since $f(x)$ has Lipschitz gradient, from \eqref{relative_inexact}, \eqref{inexact_model_1}, we can derive that after $k$  iterations (where $k$ is greater than $\sqrt{L_f /\mu}$ by a logarithmic factor $\log\left(L_f R^2/\varepsilon\right)$ with $\varepsilon$ being the desired accuracy in terms of the objective residual): 
  \begin{align}\label{rel}
    & 
    f(x_k) - f(x^*) \stackrel{\eqref{convergenceA},\eqref{deltamu}}\simeq \frac{\varepsilon}{2} + \sqrt{\frac{L_f}{\mu}}\frac{\delta_{\eqref{inexact}}^2}{\mu} \simeq  \sqrt{\frac{L_f}{\mu}}\frac{\delta_{\eqref{inexact}}^2}{\mu} \nonumber \\ 
&\hspace{0cm}
     \stackrel{\eqref{relative_inexact}, \eqref{inexact_model_1}}\simeq \sqrt{\frac{L_f}{\mu}}\frac{\alpha^2\max_{t=1,...,k}\|\nabla f(x_t)\|_2^2}{\mu}\le
     \sqrt{\frac{L_f}{\mu}}\frac{2L_f \alpha^2\max_{t=1,...,k}(f(x_k) - f(x^*))}{\mu} \nonumber \\ 
&\hspace{0cm}\lesssim
     \sqrt{\frac{L_f }{\mu}}\frac{4L_f \alpha^2\left(f(x_0) - f(x^*)\right)}{\mu}.
 \end{align}
Choosing  $\alpha = O\left( \left(\frac{\mu}{L_f}\right)^{3/4}\right)$, we guarantee that the following restart condition holds
$$
  f(x_k) - f (x^*) \leqslant \frac{1}{2}\left(f(x_0) - f (x^*)\right).
$$ 
When the restart condition holds, we restart the method. Then, after $\log\left(\Delta f/\varepsilon\right)$ restarts we can guarantee the desired $\varepsilon$-accuracy in terms of the objective residual. In ill-conditioned setting, i.e.,  when $\mu$ is small, the calculations are more involved. Yet, the main idea remains the same and replacing $\sqrt{L_f/\mu}$ with $k$ (cf. \eqref{convergenceA}) we obtain that the inequality $\alpha_k\lesssim \left(\frac{1}{k}\right)^{3/2}$ allows us to obtain the same convergence rate as in the exact gradients case.
 
 
Among many types of accelerated gradient methods, we choose to consider methods with one projection (Similar Triangles Methods (STM)), see \cite{gasnikov2018universal,cohen2018acceleration,gorbunov2019optimal,stonyakin2021inexact,dvurechensky2021first} and references therein. We choose this type of accelerated methods since: 1) it is primal-dual \cite{dvurechensky2018computational,gasnikov2018universal}; 2) it is possible to bound $\tilde{R}$ in the absence of noise \cite{gasnikov2018universal,nesterov2018lectures,stonyakin2021inexact} and when the noise is present \cite{gorbunov2018accelerated,gorbunov2019optimal}; 3) has previously been intensively investigated, see \cite{dvurechensky2021first} and references therein. 
 
 \section{Some motivation for inexact gradients}\label{motivation}
 
 In this section, we describe two, among many others, research directions where inexact gradients play an important role. We emphasise that, although the results below are not new, the way they are presented is of some value in our opinion and can be useful for the specialists in these directions.
 
 \subsection{\textbf{Gradient-free methods}}
In this subsection, we consider convex optimization problem:
\begin{equation*}
    \min_{x \in Q\subseteq \RR^n} f(x),
\end{equation*}
where $Q$ is a convex and closed set.
In some applications we do not have an access to the gradient $\nabla f(x)$ of the objective function, but we can calculate the value
\footnote{The approach we describe requires that the function values are available not only in $Q$, but also in some (depends on a particular approach) vicinity of $Q$. This problem can be solved in two different ways. The first one is ``slightly shrink the feasible set'' approach \cite{beznosikov2020gradient}. The second one is ``continuation'' of $f$ to $\RR^n$ preserving   its convexity and Lipschitz continuity \cite{risteski2016algorithms}: $f_{new}(x):= f\left(\text{proj}_Q(x)\right) + \alpha\min_{y\in Q}\|x - y\|_2$.} of
$f(x)$ with accuracy $\delta_f$ \cite{conn2009introduction}, i.e., one can evaluate $\tilde{f}(x)$ s.t.  
$$
|\tilde{f}(x) - f(x)|\leqslant \delta_f.
$$ 
\pd{An interesting question in this setting is as follows. If the accuracy $\delta_f$ of the approximation can be controlled, how should it be chosen in order to guarantee a desired accuracy $\varepsilon$ when solving problem \eqref{eq:pr_st}? A related question is what is the largest level of noise  $\delta_f$ such that the algorithm can still achieve a desired accuracy $\varepsilon$?}

In the considered setting, a number of options exists for approximating the gradient, see, e.g., \cite{berahas2021theoretical} and references therein. We consider the following examples, assuming that $f$ has $L_p$-Lipschitz $p$-th order derivatives w.r.t. the \pd{Euclidean} norm.
 \begin{itemize}
     \item \textbf{($p$-th order finite-differences).} In this case, the gradient approximation is constructed via finite differences of inexact values $\tilde{f}(x)$, which, e.g., in the case of $p=2$ lead to the following approximation to the $i$-th partial derivative
     $$
     \tilde{\nabla}_i f(x) = \frac{\tilde{f}(x + he_i) - \tilde{f}(x - he_i)}{2h}, \;\; i \in \{1,...,n\},
     $$
     where $e_i$ is the $i$-th coordinate vector and $h>0$ is a parameter. 
     For general values of $p$, we have that \eqref{inexact} holds with
     $$
     \delta = \sqrt{n}O\left(L_p h^p + \frac{\delta_f}{h} \right),
     $$
     see \cite{berahas2021theoretical}. The optimal choice of $h$ guarantees that $\delta \pd{=}  O\left(\sqrt{n}\delta_{f}^{\frac{p}{p+1}}\right)$. 
     From Section~\ref{Introduction}, we know that it is possible to solve problem \eqref{eq:pr_st} with accuracy $\varepsilon = O(\delta)$ in terms of the objective value. Hence, in order to guarantee $\varepsilon$-accuracy, we should choose $$\delta_f\sim \left(\frac{\varepsilon}{\sqrt{n}}\right)^{\frac{p+1}{p}}.$$ Unfortunately, such a simple idea does not allow one to reach the following lower bound in the class of algorithms that have sample complexity \pd{$O\left(\frac{n^{c_1}}{\varepsilon^{c_2}}\right)$}, for some $c_1,c_2 \geqslant 1$: \cite{risteski2016algorithms}  
     \begin{equation}\label{low_bound}
       \delta_f\sim \max\left\{\frac{\varepsilon^2}{\sqrt{n}},\frac{\varepsilon}{n}\right\}.  
     \end{equation}
     Note that, instead of the finite-difference approximation approach, in some applications we can use the kernel approach \cite{polyak1990optimal,bach16highly-smooth} which has recently a got renowned interest \cite{akhavan2020exploiting,novitskii2021improved}.
     \item \textbf{(Gaussian Smoothed Gradients).} In this case, the approximate gradient is formally defined as
    $$\tilde{\nabla} f(x) = \frac{1}{h}\EE{\tilde{f}(x + he)e} ,$$
    where $e\in N(0,I_n)$ is the standard Gaussian random vector. This implies that \eqref{inexact} holds with $$\delta = O\left(n^{p/2}L_p h^p + \frac{\sqrt{n}\delta_f}{h} \right),$$ see \cite{nesterov2017random,berahas2021theoretical}. The optimal choice of $h$ guarantees that $\delta \pd{=} O\left(\left(n\delta_{f}\right)^{\frac{p}{p+1}}\right)$.
    Hence, in order to guarantee $\varepsilon$-accuracy, we should choose $$\delta_f= O\left( \frac{\varepsilon^{\frac{p+1}{p}}}{n}\right).$$ 
    This bound does not match  the  lower  bound \eqref{low_bound} as well. Moreover, here (and in the next approach) we have an additional difficulty since $\tilde{\nabla} f(x)$, in general,  is not possible to evaluate exactly and only an inexact approximation is possible, for example, by the Monte Carlo approach \cite{berahas2021theoretical}, which leads to additional computational price for the better quality of approximation.
    \item \textbf{(Sphere Smoothed Gradients).}  In this case, the approximate gradient is formally defined as
    $$\tilde{\nabla} f(x) = \frac{n}{h}\EE{\tilde{f}(x + he)e} ,$$
    where $e$ is random vector with uniform distribution in the unit sphere  in $\RR^n$ with the center at $0$. This implies that \eqref{inexact} holds with $$\delta = O\left(L_p h^p + \frac{n\delta_f}{h} \right),$$ see \cite{berahas2021theoretical}. The optimal choice of $h$ guarantees $\delta = O\left( \left(n\delta_{f}\right)^{\frac{p}{p+1}}\right)$. 
    Hence, in order to guarantee $\varepsilon$-accuracy, we should choose $$\delta_f = O\left( \frac{\varepsilon^{\frac{p+1}{p}}}{n}\right).$$ 
    This bound does not match  the  lower  bound \eqref{low_bound} as well. It may seem that the this and the previous approaches are almost the same, but below we give a more accurate result for the Sphere smoothing. We are not aware of a way to obtain such a result for the Gaussian smoothing. The result is as follows \cite{devolder2014first,risteski2016algorithms}: for the Sphere smoothed gradient, we have that \eqref{inexact_model_1} holds with
    \begin{equation}\label{lower}
    \delta \simeq 2L_0h + \frac{\sqrt{n}\delta_f\tilde{R}}{h},
     \end{equation}
    where $L_0$ is the Lipschitz constant of $f$ and in \eqref{inexact_model_1} $L_f = \min\left\{L_1,\frac{7L_0^2}{h}\right\}$  when $p = 1$ and $L_f = \frac{7L_0^2}{h}$, when $p = 0$. The bound \eqref{lower} is more accurate than the previous bounds since it corresponds to the first part of the lower bound \eqref{low_bound}. Indeed, by choosing a proper $h$ in \eqref{lower} we obtain $\varepsilon\sim\delta\sim n^{1/4}\delta_{f}^{1/2}$. Hence, in order to guarantee $\varepsilon$-accuracy, we should choose $$\delta_f = O\left( \frac{\varepsilon^2}{\sqrt{n}}\right).$$
    The other part of the lower bound \eqref{low_bound}, i.e., the case when $\delta_f =O\left( \frac{\varepsilon}{n}\right)$,  is also tight, see \cite{belloni2015escaping}. Here we can also repeat the remark that the sphere smoothed gradient approximation $\tilde{\nabla} f(x)$, in general, is not available and needs to be approximated by a stochastic inexact gradient. In Section \ref{section_extinsons}, we describe an extension of our analysis of accelerated gradient method with absolute noise in the gradient to the setting of stochastic gradients.
 \end{itemize}
 
 The bound in \eqref{lower} and its consequences additionally illustrate that the inexactness and algorithms we describe in Section~\ref{Ideas} and develop below are also tight (optimal) enough. Otherwise, it would not be  possible to achieve the lower bound using the reduction of gradient-free methods to gradient methods with inexact oracle and the proposed analysis of the error accumulation for gradient-type methods.

 \subsection{\textbf{Inverse problems}}
 \label{S:inverse_problems}
 Another rather important research direction where gradients are typically available only approximately is optimization in  Hilbert spaces \cite{vasilyev2002optimization}, arising, in particular, in inverse problems theory \cite{kabanikhin2011inverse}.
 
 We start by recalling a way to calculate a derivative in a general Hilbert space. Let $J(q):=J(q,u(q))$, where $u(q)$ is the unique solution of the equation  $G(q,u)=0$. Assume that the partial $q$-derivative $G_{q}(q,u)$ of the operator $G(q,u)$ is invertible. Then, we have 
 $$G_q(q,u) + G_u(q,u)\nabla u(q) = 0, \; \text{ and } \; \nabla u(q) = - \left[G_u(q,u)\right]^{-1}G_q(q,u).$$
 Therefore,
 $$\nabla J(q)=J_q(q,u) + J_u(q,u)\nabla u(q) = J_q(q,u) - J_u(q,u)\left[G_u(q,u)\right]^{-1}G_q(q,u).$$
The same result can be obtained by considering the Lagrange functional
$$L(q,u;\psi) = J(q,u(q)) + \langle \psi, G(q,u) \rangle$$ with $$L_u(q,u;\psi)=0, \; G_q(q,u)=0, \; \text{ and } \;\nabla J(q) = L_q(q,u;\psi).$$
Indeed, by simple calculations, we can connect these two approaches by setting 
$$\psi(q,u) = - \left[G_u(q,u)^T\right]^{-1}J_u(q,u)^T.$$

Next, we demonstrate this technique on an inverse problem based on an elliptic initial-boundary-value problem.
Let $u(x,y)$ be the solution of the following problem, which we refer to as (P)
\begin{align*}
u_{xx} +u_{yy} &=0, \quad x,y\in \left( {0,1} \right), \\
u\left( {1,y} \right)&=q\left( y \right), \quad y\in \left( {0,1} \right), \\
u_x \left( {0,y} \right)&=0, \quad y\in \left( {0,1} \right), \\
u\left( {x,0} \right)&=u\left( {x,1} \right)=0, \quad x\in \left( {0,1} \right).
\end{align*}
Here we use subscripts $x,y$ to denote the corresponding partial derivatives.
The first two relations
constitute the system of equations $G(q,u) = \bar{G}\cdot(q,u)= 0$, and the last two ones constitute the feasible set $Q$.

Assume that the goal is to solve an inverse problem of estimating $q(y)\in L_2(0,1)$ by observing $b(y) = u(0,y)\in L_2(0,1)$, where $u(x,y)\in L_2\left((0,1)\times(0,1)\right)$ is the (unique) solution of (P) \cite{kabanikhin2011inverse}. We can reduce this problem to an optimization problem \cite{kabanikhin2011inverse}:
\begin{equation}\label{inverse}
\min_{q} \left\{\mathfrak{J}(q) := \min_{u:\text{ }\bar{G}\cdot(q,u)= 0, u\in Q} J(q,u):= J(u) = \int_0^1 |u(0,y) - b(y)|^2 dy\right\},
\end{equation}
which can be solved numerically since it is a convex quadratic optimization problem.
We can also directly apply Lagrange multipliers principle to \eqref{inverse}, see \cite{vasilyev2002optimization}. For that we introduce Lagrange multipliers $\psi:= \left(\psi(x,y),\lambda(y)\right)$ and write the Lagrange function:
\begin{multline*}
    L\left(q,u;\psi,\lambda\right) = J(u) + \langle \psi, \bar{G}\cdot(q,u) \rangle =
     \int_0^1 |u(0,y) - b(y)|^2 dy - \\ 
    \int_0^1\int_0^1\left(u_{xx}+u_{yy}\right)\psi(x,y)dxdy + 
    \int_0^1 \left(q(y) - u(1,y) \right)\lambda(y)dy.
\end{multline*}
To obtain a conjugate problem for $\psi$, we need to vary $L\left(q,u;\psi\right)$ in $\delta u$ satisfying $u \in Q$:
\begin{multline}\label{du}
 \delta_u L\left(q,u;\psi\right) =  2\int_0^1 \left(u(0,y) - b(y)\right)\delta u(0,y) dy -  \\
    \int_0^1\int_0^1\left(\delta u_{xx}+\delta u_{yy}\right)\psi(x,y)dxdy  -
    \int_0^1 \delta u(1,y)\lambda(y)dy,
\end{multline}
where
\begin{align*}
\delta u_x \left( {0,y} \right)&=0, \quad y\in \left( {0,1} \right),  \\
\delta  u\left( {x,0} \right)=\delta u\left( {x,1} \right)&=0,
\quad
x\in \left( {0,1} \right).
\end{align*}
Using the integration by part, from \eqref{du}, we derive
\begin{multline*}
 \delta_u L\left(q,u;\psi\right) =  \int_0^1 \left(2\left(u(0,y) - b(y)\right) - \psi_x(0,y)\right)\delta u(0,y) dy -  \\
 \int_0^1 \psi(1,y)\delta u_x(1,y) dy -  \int_0^1 \psi(x,1)\delta u_y(x,1) dx + 
  \int_0^1 \psi(x,0)\delta u_y(x,0) dy +
 \\
    \int_0^1\int_0^1\left(\psi_{xx}+  \psi_{yy}\right)\delta u(x,y)dxdy  +
    \int_0^1 \left(\psi_x(1,y) - \lambda(y)\right)\delta u(1,y)dy.
\end{multline*}

Consider now the corresponding conjugate problem, which we refer to as (D):
\begin{align*}
\psi _{xx} +\psi _{yy} &=0, \quad x,y\in \left( {0,1} \right), \\
\psi _x \left( {0,y} \right)&=2\left(u(0,y) - b(y)\right), \quad y\in \left( {0,1} \right), \\
\psi _x \left( {0,y} \right)&=2\left(u(0,y) - b(y)\right), \quad y\in \left( {0,1} \right), \\
\psi \left( {x,0} \right) &=\psi \left( {x,1} \right)=0, \quad x\in \left( {0,1} \right)
\end{align*}
and additional relation between Lagrange multipliers 
\begin{equation}\label{psi}
\lambda(y) = \psi_x(1,y),\quad y\in \left( {0,1} \right).
\end{equation}
These relations appear since $\delta_u L\left(q,u;\psi\right) = 0$,
and $\delta u(0,y),\delta u_x(1,y),\delta u(1,y) \in L_2(0,1)$;
$\delta u_y(x,1),\delta u_y(x,0) \in L_2(0,1)$;
$\delta u(x,y) \in L_2\left((0,1)\times(0,1)\right)$ are arbitrary.    

Since by \cite{rockafellar1970convex} it holds that
\[
\mathfrak{J}(q) = \min_{u: (q,u)\in (P)} J(u) = \min_{u:\text{ }\bar{G}\cdot(q,u)= 0, u\in Q} J(u) =  \min_{u\in Q}\max_{\psi\in(D)} L(q,u;\psi),
\]
from the Demyanov--Danskin's theorem \cite{rockafellar1970convex}, we have\footnote{The same result in a more simple situation (without additional constraint $u\in Q$) we considered at the beginning of this section. In that case we do not apply Demyanov--Danskin's theorem and use the inverse function theorem.}
\[
\nabla \mathfrak{J}(q) = \nabla_q \min_{u\in Q}\max_{\psi\in(D)} L(q,u;\psi) = L_q(q,u(q);\psi(q)),
\]
where $u(q)$ is the solution of (P) and $\psi(q)$ is the solution of (D), where 
\[\psi _x \left( {0,y} \right)=2\left(u(0,y) - b(y)\right),
\quad
y\in \left( {0,1} \right)
\]
and $u(0,y)$ depends on $q(y)$ via (P) and, at the same time, the pair $\left(u(q),\psi(q)\right)$ is the solution of the saddle-point problem $$\min_{u\in Q}\max_{\psi\in(D)} L(q,u;\psi).$$
Since $\delta_{\psi} L(q,u;\psi) = 0$ entails $\bar{G}\cdot(q,u)= 0$, that is from (P), if we add $u\in Q$ and $\delta_u L_(q,u;\psi) = 0$, then $u\in Q$ entails (D) as we have shown above. 
Note also that 
\[L_q(q,u(q);\psi(q))(y) = \lambda(y), \quad
y\in \left( {0,1} \right).\]
Hence, by \eqref{psi},  we have that
\[
\nabla \mathfrak{J}(q)(y) = \psi_x(1,y),\quad y\in \left( {0,1} \right).
\]
Thus we reduced the calculation of  $\nabla \mathfrak{J}(q)(y)$ to the solution of two correct initial-boundary-value problems for elliptic equation on a square, namely problems (P) and (D) \cite{kabanikhin2011inverse}.

This result can be also interpreted in a slightly different manner if we introduce a linear operator
$$
A:\quad q(y):=u(1,y)\mapsto u(0,y).
$$
Here $u(x,y)$ is the solution of problem (P). It was shown in \cite{kabanikhin2011inverse} that 
$$
A: L_2(0,1) \rightarrow  L_2(0,1).
$$
The conjugate operator is \cite{kabanikhin2011inverse}
$$
A^\ast:\quad p(y):=\psi_x(0,y)\mapsto \psi_x(1,y), \qquad 
A^\ast : L_2(0,1)\rightarrow   L_2(0,1).
$$
Here $\psi(x,y)$ is the solution of the conjugate problem (D). 
Thus, considering
\[
\mathfrak{J}(q)(y) = \|Aq - b\|_2^2, 
\]
we have
\[
\nabla \mathfrak{J}(q)(y) = A^\ast\left( 2\left(Aq - b\right)\right), 
\]
which completely corresponds to the scheme as described above:

\textbf{1.} \textit{Based on $q(y)$ we solve (P) and obtain $u(0,y) = Aq(y)$ and define $p(y) = 2\left(u(0,y) - b(y)\right)$.}

\textbf{2.} \textit{Based on $p(y)$ we solve (D) and calculate $\nabla \mathfrak{J}(q)(y) = A^\ast p(y) = \psi_x(1,y)$.}

Summarizing, the inexactness in the gradient $\nabla \mathfrak{J}(q)$ arises since we can solve (P) and (D) only numerically up to some accuracy.

The described above technique can be applied to many different inverse problems \cite{kabanikhin2011inverse} and optimal control problems \cite{vasilyev2002optimization}. Note that, for optimal control problems, in practice another strategy is widely used. Namely, instead of approximate calculation of the gradient, optimization problem is replaced by an approximate one (for example, by using finite-differences schemes). For this approximate (finite-dimensional) problem the gradient is typically available precisely \cite{evtushenko2013optimization}. Moreover, in \cite{evtushenko2013optimization} the described above Lagrangian approach is used to explain the core of automatic differentiation, where the function calculation tree is represented as a system of explicitly solvable interlocking equations.

\section{Absolute noise in the gradient} \label{section 4}
In this section, we consider problem \eqref{eq:pr_st} in the absolute noise setting (see \eqref{inexact}), i.e., we assume that the inexact gradient $\tilde{\nabla} f(x)$ satisfies uniformly in $x \in Q$  the inequality
\begin{equation}
\label{eq2.1:ref}
    \|\tilde \nabla f(x) - \nabla{f(x)} \|_2 \leqslant \delta.
\end{equation}
We underline that \(Q \) can be unbounded, for example \( \mathbb{R}^n \).
Under this assumption, we present several important relations concerning ``inexact smoothness'' and ``inexact strong convexity''. Then, we present and analyze an accelerated gradient method, study its error accumulation, and propose a stopping rule.

\subsection{Auxiliary facts}\label{sec2}
We start with some auxiliary facts and assumptions.
Let $x_{start}$ be some starting point for an algorithm and assume that there is a constant $R$ such that 
\begin{equation*}
     \|x_{start} - x^{*} \|_2 \leqslant R,
\end{equation*}
where $x^{*}$ is a solution to problem \eqref{eq:pr_st}.
If $x^{*}$ is not unique we take $x^{*}$ that is the closest to \( x_{start} \).  
We assume that the function \( f \) has Lipschitz gradient with constant \( L_f \), i.e., is $L_f$-smooth:
\begin{equation}
    \forall x, y \in Q, \; \| \nabla f(x) - \nabla f(y) \|_2 \leqslant L_f \|x - y \|_2.
    \label{eq2.2:ref}
\end{equation}
This implies the inequality
\begin{equation}
    \forall x, y \in Q, \; f(y) \leqslant f(x) + \langle \nabla f(x), y - x \rangle + \frac{L_f}{2} \|x - y \|_2^2.
    \label{eq2.3:ref}
\end{equation}
In what follows, we use \pd{the} following simple lemma.
\begin{lemma}[Fenchel inequality] \label{first lemma fehchel}
Let \( \; \space (\mathcal{E}, \langle \cdot , \cdot \rangle) \) be a Euclidean space, then
\( \forall \lambda \in \mathbb{R}_{+}, \space \forall u, v \in \mathcal{E} \),   
   \label{lemma1.1}
\begin{equation*}
     \langle u, v \rangle  \leqslant \frac{1}{2\lambda}\| u \|^2_{\mathcal{E} } + \frac{\lambda}{2} \| v \|^2_{\mathcal{E}}.
\end{equation*}
\end{lemma}
%
%
Let us introduce several constants, which will be used below in this section:
\begin{equation*}
    \begin{gathered}
        L = 2 L_f,   \qquad 
        \delta_1 = \delta,   \qquad 
        \delta_2 = \frac{\delta^2}{L}.
    \end{gathered}
\end{equation*}

From the $L_f$-smoothness assumption, we obtain the following upper bound for the objective through the inexact oracle.
\begin{claim} For all \(x, y \in Q, \) the following estimate holds:
        \label{cl2.1}
        \begin{equation*}
            f(y) \leqslant f(x) + \langle \tilde{\nabla}f(x), y - x \rangle + \frac{L}{2}\|x - y\|_2^2 + \delta_2,
        \end{equation*}
        where \(L = 2 L_f, \delta_2 = \frac{\delta^2}{2 L_f} \).
\end{claim}
\begin{proof} The proof is given by the following chain of inequalities:

    \begin{align*}
        f(y)& \leqslant f(x)
        + \langle \nabla f(x), y - x \rangle + \frac{L_f}{2} \| x - y \|_2^2 \\
        &\leqslant f(x) + \langle \tilde \nabla f(x), y - x \rangle + \frac{1}{2L_f}\|\nabla f(x) - \tilde \nabla f(x) \|_2^2 + \frac{L_f}{2} \|x - y \|_2^2 + \frac{L_f}{2} \|x - y \|_2^2 \\
        &\leqslant  f(x) + \langle \tilde \nabla f(x), y - x \rangle + \frac{L}{2} \|x - y \|_2^2 + \delta_2.
    \end{align*}
\end{proof}

We also assume that $f$ is strongly convex with parameter \( \mu \geqslant 0 \), where the case $\mu=0$ corresponds to just convexity of $f$. This means that for all \(x, y \in Q \):
\begin{equation}
    \label{eq2.4:ref}
    f(x) + \langle \nabla f(x), y - x \rangle + \frac{\mu}{2} \|x - y \|_2^2 \leqslant f(y).
\end{equation}
Based on this assumption and our assumption on the inexactness of the oracle, we can obtain two lower bounds for the objective. The first one is given by the following result.
\begin{claim} 
\label{cl2.2}
For all \(x, y \in Q \), the following estimate holds:
    \begin{equation*}
        f(x) + \langle \tilde \nabla f(x), y - x \rangle + \frac{\mu}{2} \|x - y \|_2^2 - \delta_1\|x - y \|_2 \leqslant f(y),
    \end{equation*}
    where \(\delta_1 = \delta \).
\end{claim}
\begin{proof}
    Using the Cauchy inequality and \eqref{eq2.4:ref} we obtain:
    \begin{align*}
        f(y)& \geqslant f(x) + \langle \nabla f(x), y - x \rangle + \frac{\mu}{2} \|x - y \|_2^2 \\
        &= f(x) + \langle \tilde \nabla f(x), y - x \rangle + \frac{\mu}{2} \|x - y \|_2^2
         - \langle \tilde \nabla f(x) - \nabla f(x), y - x \rangle \\ 
        & \geqslant f(x) + \langle \tilde \nabla f(x), y - x \rangle + \frac{\mu}{2} \|x - y \|_2^2 - \|\tilde \nabla f(x) - \nabla f(x) \|_2 \|x - y \|_2 \\
        & \geqslant f(x) + \langle \tilde \nabla f(x), y - x \rangle + \frac{\mu}{2} \|x - y \|_2^2 - \delta_1 \|x - y \|_2
    \end{align*}
\end{proof}
For the second estimate, we assume that $\mu \ne 0$ and introduce 
\begin{equation*}
    \delta_3 = \frac{\delta^2}{\mu}, \; \mu \not=0.
\end{equation*}
\begin{claim}\label{cl2.33} For all \(x, y \in Q \), if in \eqref{eq2.4:ref} \(\mu \not= 0 \), the following estimate holds: 
    \begin{equation*}
        f(x) + \langle \tilde \nabla f(x), y - x \rangle + \frac{\mu}{4} \|y - x \|_2^2 - \delta_3 \leqslant f(y),
    \end{equation*}
    where \(\delta_3 = \frac{\delta^2}{\mu} \).
\end{claim}
\begin{proof} Clearly,
     \begin{align*}
            f(x) & + \langle \tilde \nabla f(x), y - x \rangle + \frac{\mu}{4} \|x - y \|_2^2 - \delta_3 \\
            & = f(x) + \langle \nabla f(x), y - x \rangle  + \langle \tilde \nabla f(x) - \nabla f(x), y - x \rangle + \frac{\mu}{4} \|x - y \|_2^2 - \delta_3.
    \end{align*}
        Using Lemma~\ref{first lemma fehchel}, we obtain:
        \begin{align*}
             f(x) & + \langle \tilde \nabla f(x), y - x \rangle + \frac{\mu}{4} \|x - y \|_2^2 - \delta_3 \\
             &\leqslant f(x) + \langle \nabla f(x), y - x \rangle +  \frac{\delta^2}{\mu} + \frac{\mu}{4} \|x - y\|_2^2 + \frac{\mu}{4}\|y - x\|_2^2 - \delta_3   \leqslant f(y).
        \end{align*}
\end{proof}

To unify the derivations based on Claims \ref{cl2.2} and  \ref{cl2.33}, we use the notation $\mu_{\tau}$, $\tau \in \{1,2\}$, where $\tau = 1$ and  $\mu_1 = \mu$ correspond to the bound in Claim \ref{cl2.2} and $\tau = 2$ and  $\mu_2 = \frac{\mu}{2}$ correspond to the bound in Claim \ref{cl2.33} and the case when $\mu \ne 0$.

%
%
%
\subsection{Similar Triangles Method and its properties}\label{S:STM}
\pd{In this section, we introduce a variant of accelerated gradient method called Similar Triangles Method (STM). The design of STM is similar to that of the algorithm in \cite{gasnikov2018universal} with the main difference being that here we use inexact gradient with absolute inexactness instead of exact gradient. This change required us to  modify accordingly the analysis in order to take into account the presence of absolute inexactness in the gradient and possible unboundedness of the feasible set $Q$.}

\begin{algorithm}[H]
\caption{STM $(L, \mu_{\tau}, x_{start})$, $Q \subseteq \mathbb{R}^n$}
	\label{alg3}
\begin{algorithmic}[1]
\State 
\noindent {\bf Input:} Starting point $x_{start}$, number of steps $N$.
\State {\bf Set} $\tilde x_0 = x_{start}$, $\alpha_0 = \frac{1}{L}$, $A_0 = \frac{1}{L}$.
\State {\bf Set} $\psi_0(x) = \frac{1}{2} \|x - \tilde{x}_0\|_2^2 + \alpha_0\left(f(\tilde{x}_0) + \langle \tilde{\nabla} f(\tilde{x}_0), x - \tilde{x}_0 \rangle + \frac{\mu_{\tau}}{2} \|x - \tilde{x}_0 \|_2^2 \right)$. 
\State {\bf Set} $z_0 = \argmin_{y \in Q}{\psi_0(y)}$, $x_0 = z_0$.
\For {$k = 1 \dots N$}
        \State Find $\alpha_k$ from $(1 + \mu_{\tau} A_{k - 1})(A_{k - 1} + \alpha_k) = L \alpha_k^2$, \label{step:alpha_k}\\
        \State or equivalently $\alpha_k = \frac{1 + \mu_{\tau} A_{k - 1}}{2L} + \sqrt{\frac{(1 + \mu_{\tau} A_{k - 1})^2}{4L^2} + \frac{A_{k - 1}(1 + \mu_{\tau} A_{k - 1})}{L}}$, 
        \State $A_k = A_{k - 1} + \alpha_k,$
        \State $\tilde x_k = \frac{A_{k - 1} x_{k - 1} + \alpha_k z_{k - 1}}{A_k},$
        \State $\psi_k(x) = \psi_{k - 1}(x) + \alpha_k \left(f(\tilde{x}_k) + \langle \tilde{\nabla} f(\tilde{x}_k), x - \tilde{x}_k \rangle + \frac{\mu_{\tau}}{2} \|x - \tilde{x}_k \|_2^2 \right),$
        \State $z_k = \argmin_{y \in Q} \psi_k(y),$ 
        \State $x_k = \frac{A_{k - 1} x_{k - 1} + \alpha_k z_k}{A_k}.$
\EndFor
\State 
\noindent {\bf Output:} $x_N$.
\end{algorithmic}
\end{algorithm}
%

\begin{center}
\begin{figure}[h]
    \includegraphics[width= 0.8\linewidth]{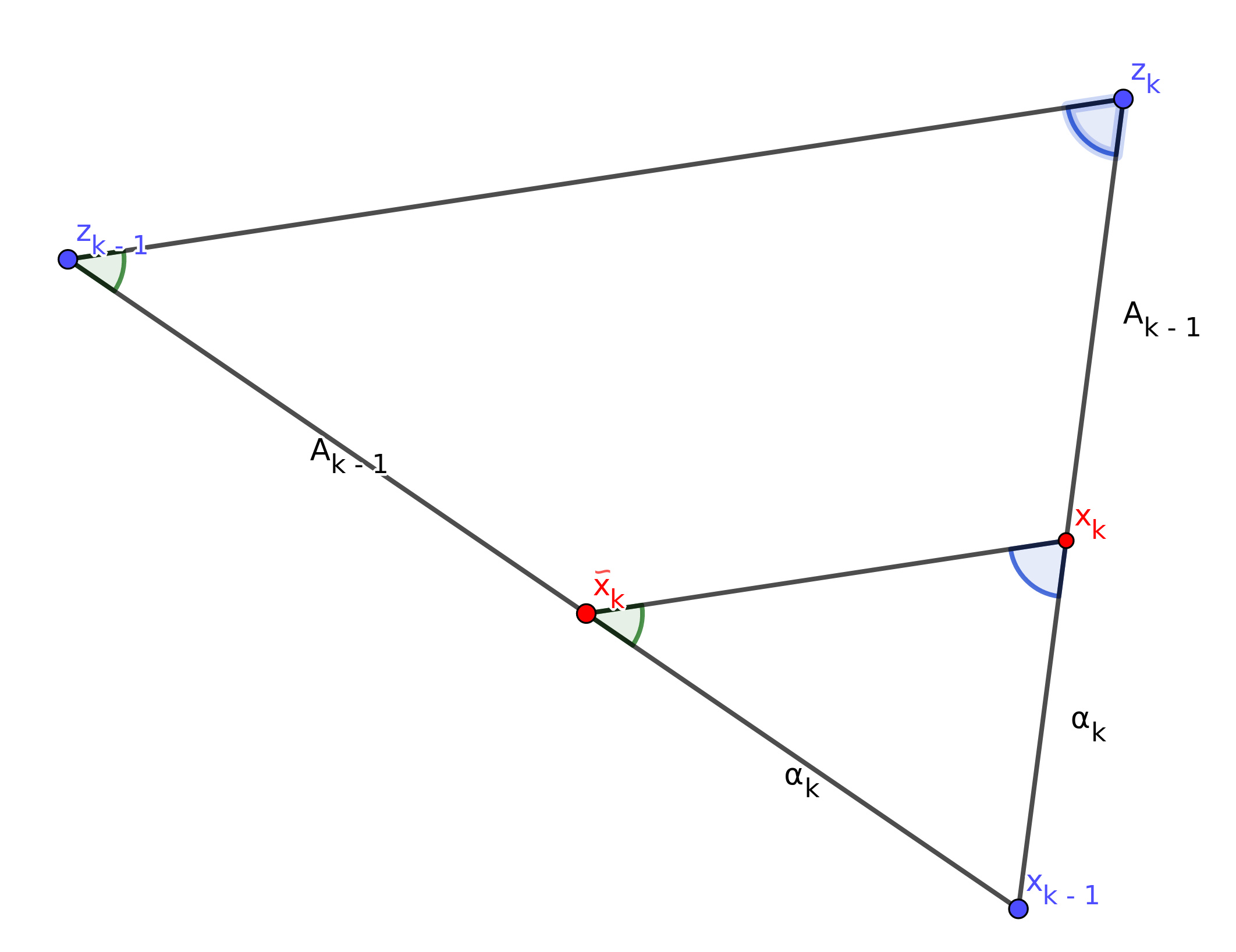}
    \caption{\ag{Geometry of Similar Triangles Method, Algorithm \ref{alg3}.}}
    \label{fig1}
\end{figure}
\end{center}
Figure~\ref{fig1} illustrates the iterates of the algorithm and justifies the name Similar Triangles Method (STM): \pd{by construction $x_k-\tilde{x}_k=\frac{\alpha_k}{A_k}(z_k-z_{k-1})$, i.e., the triangles $(z_{k-1},x_{k-1},z_k)$ and $(\tilde{x}_k,x_{k-1},x_k)$ are similar}. 
When $Q = \mathbb{R}^n$, the main step of the algorithm can be simplified to
\begin{equation*}
    z_k = z_{k - 1} - \frac{\alpha_k}{1 + \ag{A_k}\mu_{\tau}}\left(\tilde \nabla f(\tilde x_k) \ag{+} \mu_{\tau}\ag{\left(z_{k-1} - \tilde x_k\right)}\right)
\end{equation*}
using the first-order optimality condition in the definition of the point $  z_k  $.
\pd{This method is quite simple to implement, since it requires only one projection, which can be eliminated in the absence of constraints, and also has a geometric interpretation. Functions $\psi_k(x)$ contain first-order information, and are also chosen in such a way that the inequalities guaranteed by convexity or strong convexity can be used to estimate the objective from below, providing an estimating functions sequence. Moreover, since the functions $\psi_k(x)$ accumulate the first-order information from the previous iterations, the update of the variable $z_k$ can be seen as a momentum step that leads to the accelerated convergence rate. As it will be seen in Remark~\ref{Nonsmooth_structure} and Section~\ref{stoh_optim_section}, this method can be modified for composite nonsmooth optimization problems and stochastic problems.}




In the analysis, we use the following identities that easily follow from the construction of the algorithm:
\begin{equation}
\label{useful_def}
\begin{gathered}
    A_k(x_k - \tilde x_k) = \alpha_k(z_k - \tilde x_k) + A_{k - 1} (x_{k - 1} - \tilde x_k), \\
    \frac{1 + \mu_{\tau} A_{k - 1}}{2A_k} \|z_k - z_{k - 1} \|_2^2  = \frac{L}{2} \|x_k - \tilde x_k \|_2^2,  \\ 
    A_{k - 1} \|\tilde x_k - x_{k - 1} \|_2 = \alpha_k \|\tilde x_k - z_{k - 1} \|_2.
\end{gathered}
\end{equation}

The following is the main technical result which will be used later in the analysis.
\begin{lemma} For all \( k \geqslant 1 \), the following inequality holds:
    \label{LEM}
    \begin{multline*}
        \psi_k(z_k) \geqslant \psi_{k - 1}(z_{k - 1}) + \frac{1 + \mu_{\tau} A_{k - 1}}{2} \|z_k - z_{k - 1} \|_2^2
        \\
        + \alpha_k \left( f(\tilde x_k) + \langle \tilde \nabla f(\tilde x_k), z_k - \tilde x_k \rangle  + \frac{\mu_{\tau}}{2} \|z_k - \tilde x_k \|_2^2 \right).
    \end{multline*}
\end{lemma}
\begin{proof}
    By the definition of $\psi_k$, we have
    \begin{align}
        \psi_k(z_k)& = \psi_{k - 1}(z_k) + \alpha_k\left(f\left(\tilde x_k\right) + \langle \tilde \nabla f\left(\tilde x_k\right), z_k - \tilde x_k \rangle + \frac{\mu_{\tau}}{2} \|z_k - \tilde x_k \|_2^2 \right)  \notag \\
        & = \frac{1}{2} \|z_k - \tilde x_0 \|_2^2 + \displaystyle\sum_{j = 0}^{k - 1} \alpha_j\left(f\left(\tilde x_j\right) + \langle \tilde \nabla f\left(\tilde x_j\right), z_k - \tilde x_j \rangle + \frac{\mu_{\tau}}{2} \|z_k - \tilde x_j \|_2^2\right)  \notag \\
        & \hspace{1em}+ \alpha_k\left(f\left(\tilde x_k\right) + \langle \tilde \nabla f\left(\tilde x_k\right), z_k - \tilde x_k \rangle + \frac{\mu_{\tau}}{2} \|z_k - \tilde x_k \|_2^2 \right). \label{eq:LEM_proof_1}
    \end{align}
    Further, by construction, the function \(\psi_{k - 1} \) has its minimum at the point \(z_{k - 1} \), which implies, by the optimality condition,
    \begin{align}
            \langle \nabla \psi_{k - 1}(z_{k - 1}), z_k - z_{k - 1} \rangle & \geqslant 0 \notag \\ 
            \Leftrightarrow \langle z_{k - 1} - \tilde{x}_0, z_k - z_{k - 1} \rangle &\geqslant \sum_{j = 0}^{k - 1} \alpha_j \langle \tilde{\nabla} f(\tilde{x}_j) + \mu_{\tau} (z_{k - 1} - \tilde{x}_j), z_{k - 1} - z_k \rangle.
        \label{eq:LEM_proof_2}
    \end{align}
    We also have the identity
    \begin{equation}
        \frac{1}{2} \| z_k - \tilde{x}_0 \|_2^2 = \frac{1}{2} \|z_{k - 1} - z_k \|_2^2 + \frac{1}{2} \|z_{k - 1} - \tilde{x}_0 \|_2^2 
        + \langle z_{k - 1} - \tilde{x}_0, z_k - z_{k - 1} \rangle.
        \label{eq:LEM_proof_3}
    \end{equation}
    Combining the above, we have
    \begin{align*}
        \psi_k(z_k) &\stackrel{\eqref{eq:LEM_proof_1}, \eqref{eq:LEM_proof_3}}{=} \frac{1}{2} \|z_{k - 1} - \tilde x_0 \|^2_2 + \langle z_{k - 1} - \tilde x_0, z_k - z_{k - 1} \rangle + \frac{1}{2} \|z_{k - 1} - z_k \|_2^2 \\
        & \hspace{2em}+ \displaystyle\sum_{j = 0}^{k - 1} \alpha_j\left(f(\tilde x_j) + \langle \tilde \nabla f(\tilde x_j), z_k - \tilde x_j \rangle + \frac{\mu_{\tau}}{2} \|z_k - \tilde x_j \|_2^2\right) \\
        & \hspace{1em} \stackrel{\eqref{eq:LEM_proof_2}}{\geqslant} \displaystyle\sum_{j = 0}^{k - 1}\alpha_j\left( \langle \tilde \nabla f(\tilde x_j) + \mu_{\tau} (z_{k - 1} - \tilde x_j), z_{k - 1} - z_k \rangle \right) \\ 
        & \hspace{2em}+ \displaystyle\sum_{j = 0}^{k - 1} \alpha_j \left( f(\tilde x_j) + \langle \tilde \nabla f(\tilde x_j), z_k - \tilde x_j \rangle + 
        \frac{\mu_{\tau}}{2} \|z_k - \tilde x_j \|_2^2 \right) \\
        & \hspace{2em}+ \alpha_k \left(f(\tilde x_k) + \langle \tilde \nabla f(\tilde x_k), z_k - \tilde x_k \rangle + \frac{\mu_{\tau}}{2} \|z_k - \tilde x_k \|_2^2 \right) \\
        & \hspace{2em}+ \frac{1}{2} \|z_{k - 1} - \tilde x_0 \|^2_2 + \frac{1}{2} \|z_{k - 1} - z_k \|_2^2. 
    \end{align*}
Applying the identity
    \begin{equation*}
         \langle z_{k - 1} - \tilde x_j, z_{k - 1} - z_k \rangle = 
         \frac{1}{2} \|z_{k - 1} - \tilde{x}_j \|_2^2 + \frac{1}{2} \|z_k - z_{k - 1} \|_2^2 - \frac{1}{2} \|z_k - \tilde{x}_j \|_2^2, 
    \end{equation*}
 and the definition of the sequence $\lbrace A_k \rbrace$, we finally get
    \begin{align*}
        \psi_k(z_k)& \geqslant \frac{1}{2} \|z_{k - 1} - \tilde x_0 \|^2_2 + \frac{1 + \mu_{\tau} A_{k - 1}}{2} \|z_{k - 1} - z_k \|_2^2 
        \\
        & \hspace{1em}+ \displaystyle\sum_{j = 0}^{k - 1} \alpha_j\left(f(\tilde x_j) + \langle \tilde \nabla f(\tilde x_j), z_{k - 1} - \tilde x_j \rangle + \frac{\mu_{\tau}}{2} \|z_{k - 1} - \tilde x_j \|_2^2 \right) \\
        & \hspace{1em}+ \alpha_k \left(f(\tilde x_k) + \langle \tilde \nabla f(\tilde x_k), z_k - \tilde x_k \rangle
        + \frac{\mu_{\tau}}{2} \| z_k - \tilde x_{k} \|_2^2 \right) 
        \\
        &= \psi_{k - 1}(z_{k - 1}) + \frac{1 + \mu_{\tau} A_{k - 1}}{2} \|z_k - z_{k - 1} \|_2^2
        \\
        & \hspace{1em}+ \alpha_k \left(f(\tilde x_k) + \langle \tilde \nabla f(\tilde x_k), z_k - \tilde x_k \rangle  + \frac{\mu_{\tau}}{2} \|z_k - \tilde x_k \|_2^2 \right)
    \end{align*}
\end{proof}
\begin{remark} 
In the case when \(\mu = 0\), we obtain the following particular case of the result of Lemma~\ref{LEM}:
    \begin{equation*}
        \begin{gathered}
        \psi_k(z_k) = \psi_{k - 1}(z_k) + \alpha_k \left( f(\tilde x_k) + \langle \tilde \nabla{f}(\tilde x_k), z_k - \tilde x_k \rangle \right) \\ 
        \Rightarrow \psi_k(z_k) \geqslant \psi_{k - 1}(z_{k - 1}) + \frac{1}{2} \|z_k - z_{k - 1}\|_2^2 + \alpha_k \left( f(\tilde x_k) + \langle \tilde \nabla{f}(\tilde x_k), z_k - \tilde x_k \rangle \right).
        \end{gathered}
    \end{equation*}
\end{remark}

We finish this subsection by a series of technical results that estimate the growth of the sequence $\lbrace A_k \rbrace$ and related sequences.
\begin{claim} If $\mu \not= 0$, then for all $ k \geqslant 1$ the following inequality holds:
    \label{cl3.1}
    \begin{equation*}
        A_k \geqslant A_{k - 1} \lambda_{\mu_{\tau}, L},
    \end{equation*}
    where
    $$\theta_{\mu_{\tau}, L} = \frac{\mu_{\tau}}{L}, \; \lambda_{\mu_{\tau}, L} =  \left(1 + \frac{1}{2}\theta_{\mu_{\tau}, L} + \sqrt{\theta_{\mu_{\tau}, L}} \right). $$
\end{claim}

\begin{proof}
    Using the relation between $\alpha_k, A_k, A_{k - 1}$: 
    \begin{equation*}
        A_k(1 + \mu_{\tau} A_{k - 1}) = L \alpha_k^2,
    \end{equation*}
we obtain a quadratic equation for $A_k$:
\begin{equation*}
    \begin{gathered}
        A_k(1 + \mu_{\tau}A_{k - 1}) = L(A_k - A_{k - 1})^2, \\
        A_k(1 + \mu_{\tau}A_{k - 1}) = LA_k^2 - 2LA_{k - 1}A_k + LA_{k - 1}^2, \\
        LA_k^2 - A_k \left( 1 + \mu_{\tau} A_{k - 1} + 2L A_{k - 1} \right) = 0
        + L A_{k - 1}^2.
    \end{gathered}
\end{equation*}
Solving this equation, we get
\begin{equation*}
    \begin{gathered}
        A_k \geqslant A_{k - 1} \left(1 + \frac{\mu_{\tau} + 1}{L} + \sqrt{\frac{\mu_{\tau}}{L}} \right) \geqslant A_{k - 1} \left(1 + \frac{\mu_{\tau}}{2L} + \sqrt{\frac{\mu_{\tau}}{L}} \right).
    \end{gathered}
\end{equation*}
\end{proof}
Using that, for $ x < 1$, $1 + x \geqslant e^{\frac{x}{2}}$, we obtain the following result.
\begin{corollary}
\label{lower bound lambda}
\begin{equation*}
    \begin{gathered}
        \lambda_{\mu_{\tau}, L} = \left(1 + \frac{\mu_{\tau}}{2L} +  \sqrt{\theta_{\mu_{\tau}, L}} \right) \geqslant \left(1 + \sqrt{\theta_{\mu_{\tau}, L}} \right) \geqslant e^{\frac{1}{2}\sqrt{\theta_{\mu_{\tau}, L}}}.
    \end{gathered}
\end{equation*}
\end{corollary} 
\begin{claim} If $\mu \not= 0$, then for all $ k \geqslant 1$ the following inequality holds:
    \label{cl3.2}
    \begin{equation*}
        \frac{1}{A_k} \displaystyle\sum_{j = 0}^{k}A_j \leqslant 1 + \sqrt{\frac{L}{\mu_{\tau}}} .
    \end{equation*}
\end{claim}
\begin{proof}
    %
    Using the previous claim, we get $A_k \geqslant A_{k - j} \lambda_{\mu_{\tau}, L} ^j$, and, hence, $\frac{A_{k - j}}{A_k} \leqslant \lambda_{\mu_{\tau}, L}^{-j}$.
    %
    This gives
    \begin{equation*}
       \frac{1}{A_k}\displaystyle\sum_{j = 0}^{k}A_j \leqslant \displaystyle\sum_{j = 0}^{k} \lambda_{\mu_{\tau}, L}^{-j} = \frac{\lambda_{\mu_{\tau}, L}^{k + 1} - 1}{\lambda_{\mu_{\tau}, L}^{k + 1} - \lambda_{\mu_{\tau}, L}^k} \leqslant \frac{\lambda_{\mu_{\tau}, L}}{\lambda_{\mu_{\tau}, L} - 1} \leqslant 1 + \sqrt{\frac{L}{\mu_{\tau}}}.
    \end{equation*}
\end{proof}
\begin{claim} If \( \mu = 0 \), then for all $ k \geqslant 1$, 
\label{cl3.3}
\begin{equation*}
    A_k \geqslant \frac{(k + 1)^2}{4L}.
\end{equation*}
\end{claim}
\begin{proof}
    
    If \(\mu = 0 \), then \(A_k = L \alpha_k^2 \), and, solving the quadratic equation, we get 
    \begin{equation*}
        \alpha_k = \frac{1 + \sqrt{1 + 4L^2\alpha_{k - 1}^2}}{2L} \geqslant \frac{1+2L\alpha_{k - 1}}{2L} = \frac{1}{2L} + \alpha_{k - 1}.
    \end{equation*}
    Then, by induction, it is easy to see that
    \begin{equation*}
       \alpha_k \geqslant \frac{k + 1}{2L} \Rightarrow A_k = L a_k^2 \geqslant \frac{(k + 1)^2}{2L}. 
    \end{equation*}
    
\end{proof}
\begin{claim} If \(\mu = 0 \), then for all $ k \geqslant 1$
\label{cl3.4}
    \begin{equation*}
        \frac{1}{A_k} \displaystyle\sum_{j = 0}^{k}A_j\leqslant k + 1.
    \end{equation*}
\end{claim}
\begin{proof} The proof follows from the simple calculations since $\lbrace A_k \rbrace$ is non-decreasing:
    \begin{equation*} 
        \frac{1}{A_k} \displaystyle\sum_{j = 0}^{k}A_j \leqslant
        \frac{1}{A_k} (k + 1)A_k = k + 1.
    \end{equation*}
\end{proof}

\subsection{Convergence rates under the absolute inexactness}
In this section we obtain main convergence rate results for Algorithm \ref{alg3}.
We will use the following sequence
\begin{equation}
    \label{equationR}
    \tilde R_k = \max\limits_{0 \leqslant j \leqslant k} \lbrace \|z_j - x^*\|_2, \|x_j - x^*\|_2, \|\tilde x_j - x^*\|_2 \rbrace.
\end{equation}

\pd{\begin{proposition} 
\label{thr4.1}
The sequences $\lbrace x_k \rbrace$, $\lbrace \tilde{x}_k \rbrace$, $\lbrace z_k \rbrace$ generated by Algorithm $\ref{alg3}$ satisfy for all $k \geqslant 0$ the inequality
    \begin{equation*}
        A_k f(x_k) \leqslant \psi_k(z_k) + \delta_2 \displaystyle\sum_{j=0}^{k} A_j + \delta_1 \displaystyle\sum_{j=1}^{k} \alpha_j \|\tilde{x}_j-z_{j-1}\|_2.  
    \end{equation*}
\end{proposition}
}
\begin{proof} We prove the result by induction. The induction basis for $k=0$ follows from the facts that $A_0=\alpha_0 = \frac{1}{L}$ and 
\begin{equation*}
        \psi_0(x) = \alpha_0 \left( f(\tilde{x}_0) + \langle \tilde{\nabla} f(\tilde{x}_0), x - \tilde{x}_0 \rangle + \frac{\mu_{\tau}}{2}
        \|x - \tilde{x}_0 \|_2^2 \right) + \frac{1}{2} \|x - \tilde{x}_0 \|_2^2,
\end{equation*}
which imply, by Claim \ref{cl2.1}, and since $x_0 = z_0$, that
\begin{multline*}
    f(x_0)
    \leqslant f(\tilde x_0) + \langle \tilde \nabla f(\tilde x_0), x_0 - \tilde x_0 \rangle + \frac{L}{2}\|x_0 - \tilde x_0 \|_2^2 + \delta_2
    \\
    = L\psi_0(z_0) - \frac{\mu_{\tau}}{2} \|z_0 - \tilde x_0 \|_2^2 + \delta_2  \leqslant L\psi_0(z_0) + \delta_2.
\end{multline*}
To make the induction step, we start from the following corollary of Claim \ref{cl2.1}:
\begin{multline*}
    A_k f(x_k) - A_{k-1} \delta_1 \|x_{k - 1} - \tilde x_k \|_2  \\
    \leqslant A_k \left( f(\tilde x_{k}) + \langle \tilde \nabla f(\tilde x_{k}), x_k - \tilde x_k \rangle + \frac{L}{2} \|x_k - \tilde x_k \|_2^2 + \delta_2 \right) - A_{k - 1} \delta_1 \|x_{k - 1} - \tilde x_k \|_2.
\end{multline*}
Using equations $\eqref{useful_def}$, this gives
\begin{align*}
    A_k f(x_k) & - A_{k-1} \delta_1 \|x_{k - 1} - \tilde x_k \|_2 \\
    & \leqslant A_{k - 1} \left( f(\tilde x_k) + \langle \tilde \nabla f(\tilde x_k), x_{k - 1} - \tilde x_k \rangle \right) \\
   & \hspace{1em}+ \alpha_k \left( f(\tilde x_k) + \langle \tilde \nabla f(\tilde x_k), z_k - \tilde x_k \rangle \right) \\
   & \hspace{1em}+ \frac{(1 + \mu_{\tau} A_{k - 1})}{ 2}\|z_k - z_{k - 1} \|_2^2  + A_k \delta_2 - A_{k - 1}\delta_1\|x_{k - 1} - \tilde x_k \|_2 \\
   &\leqslant A_{k - 1}f(x_{k - 1}) + \alpha_k(f(\tilde x_k) + \langle \tilde \nabla f(\tilde x_k), z_k - \tilde x_k \rangle ) \\
   & \hspace{1em} +  \frac{1 + \mu_{\tau} A_{k - 1}}{2}\|z_k - z_{k - 1} \|_2^2 + A_k \delta_2.
\end{align*}
By the induction hypothesis and since $\frac{\mu_{\tau}}{2}\|z_k-\tilde{x}_k\|_2^2 \geqslant 0$, we further obtain
\begin{multline*}
    A_k f(x_k) - A_{k-1} \delta_1 \|x_{k - 1} - \tilde x_k \|_2 \leqslant \psi_{k - 1}(z_{k - 1}) + \delta_2 \displaystyle\sum_{j=0}^{k - 1} A_j + \delta_1 \displaystyle\sum_{j=1}^{k-1} \alpha_j \|\tilde{x}_j-z_{j-1}\|_2 \\ 
   + \alpha_k \left( f(\tilde x_k) + \langle \tilde \nabla f(\tilde x_k), z_k - \tilde x_k \rangle +\frac{\mu_{\tau}}{2}\|z_k-\tilde{x}_k\|_2^2 \right) + \frac{1 + \mu_{\tau} A_{k - 1}}{2} \|z_k - z_{k - 1} \|_2^2  + A_k \delta_2.
\end{multline*}
Using Lemma~$\ref{LEM}$, we get
\begin{align*}
    A_k f(x_k) &\leqslant A_{k - 1}\delta_1\|x_{k - 1} - \tilde x_k \|_2 + \psi_k(z_k) + \delta_2 \displaystyle\sum_{j=0}^{k} A_j + 
    \delta_1 \displaystyle\sum_{j=1}^{k-1} \alpha_j \|\tilde{x}_j-z_{j-1}\|_2
    \\ &\stackrel{\eqref{useful_def}}{=} \psi_k(z_k) + \delta_2 \displaystyle\sum_{j=0}^{k} A_j + 
    \delta_1 \displaystyle\sum_{j=1}^{k-1} \alpha_j \|\tilde{x}_j-z_{j-1}\|_2  + \alpha_k \delta_1 \|\tilde 
    x_k - z_{k - 1} \|_2,
\end{align*}
which finishes the induction step and the proof.
\end{proof}
Using the definition of $\lbrace \tilde{R}_k \rbrace$ and $\lbrace A_k \rbrace$, we obtain the following simple corollary of the above \pd{proposition}:
\begin{align}
    A_k f(x_k)& \leqslant \psi_k(z_k) + \delta_2 \displaystyle\sum_{j=0}^{k} A_j + \delta_1 \displaystyle\sum_{j=1}^{k} \alpha_j \|\tilde{x}_j-z_{j-1}\|_2 \notag
    \\
    &\leqslant \psi_k(z_k) + \delta_2 \displaystyle\sum_{j=0}^{k} A_j + \delta_1 \displaystyle\sum_{j=1}^{k} \alpha_j (\|z_{k - 1} - x^* \|_2 + \|\tilde x_k - x^* \|_2) \notag
    \\
    &\leqslant \psi_k(z_k) + \delta_2 \displaystyle\sum_{j=0}^{k} A_j + 2 \tilde R_k \delta_1 A_k.
\label{eq:important_bound_abs_inex}
\end{align}

We note that the above estimates hold both in the case of \(\mu \not= 0\) and in the case of \(\mu = 0\).

%
The proof of the following result repeats verbatim the proof of \pd{Proposition}~\ref{thr4.1}, except for Claim~\ref{cl2.2} being replaced by Claim~\ref{cl2.33}.
\pd{
\begin{proposition}
If $ \mu \not= 0 $, the sequences $\lbrace x_k \rbrace$, $\lbrace \tilde{x}_k \rbrace$, $\lbrace z_k \rbrace$ generated by Algorithm $\ref{alg3}$ satisfy for all $k \geqslant 0$ the inequality
    \label{th4.2}
    \begin{equation*}
        A_k f(x_k) \leqslant \psi_k(z_k) + \delta_2 \displaystyle\sum_{j=0}^{k} A_j + \delta_3 \displaystyle\sum_{j=0}^{k - 1} A_j. 
    \end{equation*}
\end{proposition} 
}
\pd{
\begin{proposition} Assume that the oracle error $\delta$ satisfies \( \delta = 0 \) and that $ \|\tilde{x}_0 - x^* \|_2 \leqslant R$ for some $R$.
Then, for all $k\geqslant 1$, $\tilde R_k \leqslant R$. 
\label{th4.3}
\end{proposition}
}
\begin{proof}
    We first prove that, for all $k\geqslant 0$, $\|z_k-x^*\|_2 \leqslant R$. Let us fix $k\geqslant 0$. 
    By \pd{Proposition}~\ref{thr4.1}, we have \( \; A_k f(x_k) \leqslant \psi_k(z_k) \). Further, $\psi_k(x)$ is strongly convex with the constant at least 1.  At the same time, by the strong convexity of $f$, we have
    \begin{equation*}
        f(\tilde x_j) + \langle \nabla f(\tilde x_k), x^* - \tilde x_k \rangle  + \frac{\mu_{\tau}}{2} \|x^* - \tilde x_k \|_2^2
        \leqslant f(x^*) \leqslant f(x_k).
    \end{equation*}
  Using these three facts and the definition of $\psi_k(x)$, we obtain:
    \begin{align*}
        \frac{1}{2}\|z_k - x^* \|_2^2& = \frac{1}{2}\|z_k - x^* \|_2^2 + A_k f(x_k) - A_k f(x_k)
        \\
        &\leqslant \psi_k (z_k) + \frac{1}{2}\|z_k - x^* \|_2^2 - A_k f(x_k) \\
        &\leqslant \psi_k(x^*) - A_k f(x_k) \\
        &\leqslant \displaystyle\sum_{j = 0}^{k} \alpha_j\left(f(\tilde x_j) + \langle \nabla f(\tilde x_k), x^* - \tilde x_k \rangle + \frac{\mu_{\tau}}{2} \|x^* - \tilde x_k \|_2^2 \right) \\
        & \hspace{1em}+ \frac{1}{2} \|x^* - \tilde x_0 \|_2^2 - A_k f(x_k) \\
        & \leqslant \sum\limits_{j=0}^{k} \alpha_j f(x^*) + \frac{1}{2}R^2 - A_k f(x_k)
        \\
        &= A_k \left(f(x^*) - f(x_k) \right) + \frac{1}{2}R^2  \leqslant \frac{1}{2}R^2.
    \end{align*}

    %
    For the remaining two sequences, $\lbrace \tilde{x}_k \rbrace$ and $\lbrace x_k \rbrace$ the proof is organized by induction. Clearly, $\|\tilde{x}_0-x^*\|\leqslant R$. Since, by construction, $x_0=z_0$, we have $\|x_0-x^*\|\leqslant R$. Then, by construction of the algorithm and the induction hypothesis, we have
    \begin{align*}
         \|x_k - x^* \|_2 & = \left\|\frac{A_{k - 1}}{A_k}(x_{k - 1} - x^*) + \frac{\alpha_k}{A_k}(z_k - x^*) \right\|_2
         \\
         &\leqslant \frac{A_{k - 1}}{A_k} \|x_{k - 1} - x^* \|_2 + \frac{\alpha_k}{A_k} \|z_k - x^* \|_2 \leqslant R.
    \end{align*}
    In the same way, we obtain $\|\tilde{x}_k-x^*\|\leqslant R$ using the definition
    \begin{equation*}
        \tilde{x}_k = \frac{\alpha_k}{A_k}z_{k - 1} + \frac{A_{k - 1}}{A_k} x_{k - 1}.
    \end{equation*}
\end{proof}
\pd{Using the above results, we obtain the following convergence rate result for the STM algorithm.}
\begin{theorem}[Main Theorem]
\label{Th:STM_main}
Let $ \|\tilde{x}_0 - x^* \|_2 \leqslant R$ for some $R$.
If $ \mu \not= 0 $, the sequences \pd{$\lbrace x_k \rbrace$, $\lbrace \tilde{x}_k \rbrace$, $\lbrace z_k \rbrace$} generated by Algorithm $\ref{alg3}$ satisfy for all $\pd{N} \geqslant 0$ the inequalities
    \label{th4conv}
    \begin{equation*}
        \begin{gathered}
            f(x_N) - f(x^*) \leqslant L R^2 \exp\left({{-\frac{1}{2}\sqrt{\frac{\mu_1}{L}}N}}\right) + \left(1 + \sqrt{{\frac{L}{\mu_1}}}\right)\delta_2 + 3 \tilde R_N \delta_1, \\ 
            f(x_N) - f(x^*) \leqslant L R^2 \exp{\left(-\frac{1}{2}\sqrt{\frac{\mu_2}{L}}N \right)} + \left(1 + \sqrt{{\frac{L}{\mu_2}}} \right)\delta_2 + \left(1 + \sqrt{\frac{L}{\mu_2}} \right)\delta_3. \\
        \end{gathered}
    \end{equation*}
If $ \mu = 0 $, the sequences $\lbrace x_k \rbrace$, $\lbrace \tilde{x}_k \rbrace$, $\lbrace z_k \rbrace$ generated by Algorithm $\ref{alg3}$ satisfy for all $\pd{N} \geqslant 0$ the inequality  
\begin{equation}
\label{eq:STM_convergence_conv}
        f(x_N) - f(x^*) \leqslant \frac{4LR^2}{N^2} + 3 \tilde R_N \delta_1 + N \delta_2,
\end{equation}
where the sequence $\lbrace \tilde R_k \rbrace$ is defined in \eqref{equationR}.
\end{theorem}
\begin{proof}
    The proofs of the first and second inequalities are nearly the same with the only difference that the proof of the first inequality is based on \pd{Proposition}~\ref{thr4.1} and Claim~\ref{cl2.2}, whereas the proof of the second inequality is based on \pd{Proposition}~\ref{th4.2} and Claim~\ref{cl2.33}. Thus, we give only the proof of the first inequality.
    From \eqref{eq:important_bound_abs_inex}, by the definition of $\lbrace z_N \rbrace$ and $\lbrace \psi_N(\cdot) \rbrace$, and Claim~\ref{cl2.2}, we have
    \begin{align*}
        A_N f(x_N)& \leqslant \psi_N(z_N) + \delta_2 \displaystyle\sum_{j=0}^{N} A_j + 2 \tilde R_N \delta_1 A_N \\
        &\leqslant \frac{1}{2} \|x^* - \tilde x_0 \|_2^2
        + \delta_2 \displaystyle\sum_{j=0}^{N} A_j \\
        & \hspace{1em}+ 2 \tilde R_N \delta_1 A_N + \displaystyle\sum_{j=0}^N \alpha_k(f(\tilde x_j) + \langle \tilde \nabla f(\tilde x_j), x^* - \tilde x_j \rangle + \frac{\mu_1}{2} \|x^* - x_i \|_2^2 ) \\
        & \leqslant \delta_2 \displaystyle\sum_{j=0}^{N} A_j + 2 \tilde R_N \delta_1 A_N + \displaystyle\sum_{j = 0}^{N} \alpha_k( \tilde R_j \delta_1 + f(x^*) ) + \frac{1}{2}R^2 \\
        &= \delta_2 \displaystyle\sum_{j=0}^{N} A_j + 3 \tilde R_N \delta_1 A_N + A_N f(x^*) + \frac{1}{2}R^2 \\
        &\Longleftrightarrow
        f(x_N) - f(x^*) \leqslant \frac{R^2}{2A_N} + \delta_2 \frac{1}{A_N}\sum\limits_{j=0}^{N}A_j + 3 \tilde R_N \delta_1.
    \end{align*}
    %
    Using Claim~\ref{cl3.1} with Corollary~\ref{lower bound lambda} and Claim~\ref{cl3.2} we get:
    \begin{equation*}
        f(x_N) - f(x^*) \leqslant L R^2 \exp{\left(-\frac{1}{2}\sqrt{\frac{\mu_1}{L}}N \right)} + \left(1 + \sqrt{{\frac{L}{\mu_1}}}\right)\delta_2 + 3 \tilde R_N \delta_1.
    \end{equation*}
    
    Repeating the same steps and using Claim~\ref{cl3.3} and \ref{cl3.4}, we prove the third inequality.
\end{proof}
 
\pd{Commenting on the results obtained in Theorem \ref{Th:STM_main}, we can conclude that in the case of strong convexity and the presence of absolute noise, STM converges in  terms of the objective value up to some limiting accuracy. Namely, the convergence rate bound is the sum of the convergence rate of the optimal method for the class of strongly convex and Lipschitz-smooth problems and the term characterizing the limiting error caused by the noise is
\begin{equation*}
    \left(1 + \sqrt{{\frac{L}{\mu_2}}} \right)\left(\delta_2 + \delta_3 \right).
\end{equation*}
In the case when $\mu = 0$, we obtain a weaker convergence rate statement since in the estimate we see a linear accumulation of the noise in the term $N \delta_2$, as well as in the term $ 3\tilde R_N \delta_1$ (Note that \pd{Proposition}~\ref{th4.3} gives the estimate for $\tilde R_N$ in the absence of noise). This motivates us to use the regularization technique to make a reduction of the convex case to the strongly convex case, which is considered in the next Remark~\ref{remark4.2}. Another way to deal with the noise accumulation is to introduce a stopping rule, which is done below in Section~\ref{subsection_stop_rule}.}

%
\begin{remark}
    \label{remark4.2}
We can make a reduction of the setting when $\mu =0$ to the setting when $\mu \ne 0$. Indeed, suppose that $\mu = 0$ and consider the following regularized problem: 
\begin{equation*}
    \begin{gathered}
        \min\limits_{x \in Q} \left\{f^{\mu}(x) := f(x) + \frac{\mu}{2} \|x - \tilde x_0 \|_2^2 \right\}.
    \end{gathered}
\end{equation*}
Then, we have
\begin{equation*}
    \begin{gathered}
        \nabla f^{\mu}(x) = \nabla f(x) + \mu(x - \tilde x_0), \\
        \tilde \nabla f^{\mu}(x) = \tilde \nabla f(x) + \mu(x - \tilde x_0), \\
        \|\tilde \nabla f^{\mu}(x) - \nabla f^{\mu}(x) \|_2 = \|\tilde \nabla f(x) - \nabla f(x) \|_2 \leqslant \delta.
    \end{gathered}
\end{equation*}
Clearly, $  f^{\mu}(x)$ has Lipschitz gradient. Indeed, \(\forall x, y \in Q \):
\begin{align*}
    \|\nabla f^{\mu} (x) - \nabla f^{\mu}(y) \|_2& = \|(\nabla f(x) - \nabla f(y))  + \mu (x - y) \|_2 \\
     &\leqslant  \|\nabla f(x) - \nabla f(y) \|_2 + \mu \|x - y \|_2 \\
     &\leqslant L_f \|x - y \|_2 + \mu \|x - y \|_2 \leqslant (L_f + \mu) \|x - y \|_2.
\end{align*}
Since $\mu \leqslant L$, we have that $  f^{\mu}(x)$ is $L^\mu$-smooth with \(L^\mu = 4L_f = 2L \). Moreover, $f^{\mu}(x)$ is strongly convex and we can apply the derivations corresponding to the case $\tau=2$. Using Theorem~\ref{th4conv}, and setting $x_{\mu}^* = \argmin_{x \in Q} f^{\mu}(x)$, $R_{\mu}$ s.t. $ \|x_{\mu}^* - \tilde x_0 \|_2 \leqslant R_{\mu}$, we obtain the following inequalities
\begin{equation*}
\begin{gathered}
    f^{\mu}(x_k) - f^{\mu}(x_{\mu}^*) \leqslant 2L R_{\mu}^2 \exp\left({-\frac{1}{2}\sqrt{\frac{\mu}{4L}}k} \right) + \left( 1 + \sqrt{\frac{4L}{\mu}} \right) \left( \frac{1}{2L} + \frac{1}{\mu} \right) \delta^2, \\
    f^{\mu}(x_{\mu}^*) \leqslant f(x^*) + \frac{\mu}{2}R^2.
\end{gathered}
\end{equation*}
Translating this to the original objective $f$, we obtain
\begin{align*}
    f(x_k) - f(x^*)& \leqslant f^{\mu}(x_k) - f(x^*)
    \\
    &\leqslant f^{\mu}(x_k) - f(x_{\mu}^*) + \frac{\mu}{2}R^2
    \\
    &\leqslant 2L R_{\mu}^2 \exp\left({-\frac{1}{2}\sqrt{\frac{\mu}{4L}}k} \right) + \left( 1 + \sqrt{\frac{4L}{\mu}} \right) \left( \frac{1}{2L} + \frac{1}{\mu} \right) \delta^2 + \frac{\mu}{2}R^2.
\end{align*}
By the strong convexity of the function $f^{\mu}$, we get:
\begin{equation*}
    \begin{gathered}
        f(x^*) + \frac{\mu}{2}R_{\mu}^2 \leqslant f(x_{\mu}^*) + \frac{\mu}{2}R_{\mu}^2 = f^{\mu}(x_{\mu}^*) \leqslant f^{\mu}(x^*) = f(x^*) + \frac{\mu}{2}R^2 \Rightarrow \\
        R_{\mu} \leqslant R.
    \end{gathered}
\end{equation*}
Finally, we get the convergence rate as follows:
\begin{multline*}
    \label{conv_for_reg}
    f(x_k) - f(x^*) \leqslant 2L R^2 \exp\left({-\frac{1}{2}\sqrt{\frac{\mu}{2L}}k} \right) + \left( 1 + \sqrt{\frac{4L}{\mu}} \right) \left( \frac{1}{2L} + \frac{1}{\mu} \right) \delta^2 + \frac{\mu}{2}R^2. 
\end{multline*}
To obtain an error $\varepsilon$ in the r.h.s., we choose $\mu = \frac{2}{3} \frac{\varepsilon}{R^2}$.
\end{remark}
%
%
%
\subsection{Stopping rule under the absolute inexactness}
\label{subsection_stop_rule}
In this subsection, we consider the setting with $\tau = 1$ and $\mu = 0$.
In this case, a possible drawback of the convergence rate obtained in Theorem~\ref{th4conv}
    \begin{equation*}
        f(x_N) - f(x^*) \leqslant \frac{4LR^2}{N^2} + 3 \tilde{R}_N \delta_1 + N \delta_2
    \end{equation*}
    can be that the sequence $\lbrace \tilde{R}_N \rbrace$ may increase as $N$ increases. 
    To overcome this, we formulate a certain condition (stopping rule) and prove that if it is satisfied at iteration $N$, the algorithm solves problem \eqref{eq:pr_st} with certain accuracy, and if it is not satisfied at iteration $N$, then $\tilde{R}_N \leqslant R$. Moreover, we estimate the maximum number of iterations to satisfy this condition.
    \begin{theorem}
        \label{stopping_rule}
        Consider the setting $\tau = 1$ and $\mu = 0$ and assume that for some $R$, $\|\tilde{x}_0-x^*\|_2 \leqslant R$. Let $\varepsilon >0$ be the desired solution accuracy. Let $N$ be the first iteration such that 
        \begin{equation}
            \label{eq:stopping_rule_p}
            f(x_N)-f(x^*) \leqslant \frac{\delta_2}{A_N}\sum_{j=0}^{N}A_j + 3R \delta_1 + \varepsilon.
        \end{equation}
        Then, for all $k \in \{0,\ldots,N-1\}$, we have that $\tilde{R}_{k} \leqslant R$. Moreover,
        \begin{equation}
            \label{eq:stopping_rule_N_bound}
            N \leqslant N_{\max}:= \left\lceil \sqrt{\frac{2LR^2}{\varepsilon}}\right\rceil .
        \end{equation}
    \end{theorem}
    \begin{proof}
    Fixing any $k \geqslant 0$, applying \pd{Proposition}~\ref{thr4.1}, the fact that $1$-strongly convex function $\psi_k (\cdot)$ attains its minimum at the point $  z_k$, the definition of this function, and Claim \ref{cl2.2}, we obtain
    \begin{align}
        \frac{1}{2} \|z_k - x^* \|_2^2  + A_k f(x_k)& \notag  
        \leqslant \frac{1}{2} \|z_k - x^* \|_2 + \psi_k(z_k) + \delta_2 \displaystyle \sum_{j = 0}^{k} A_j +  \delta_1 \displaystyle \sum_{j = 1}^{k} \alpha_j \|\tilde x_j - z_{j - 1} \|_2 \notag\\
        &\leqslant \psi_k(x^*) +\delta_2 \displaystyle \sum_{j = 0}^k A_j +  \delta_1 \displaystyle \sum_{j = 1}^{k} \alpha_j \|\tilde x_j - z_{j - 1} \|_2 \notag\\
        & = \frac{1}{2}\|\tilde{x}_0-x^*\|_2^2 + \sum_{j=0}^k \alpha_j (f(\tilde{x}_j) + \langle \tilde \nabla f\left(\tilde x_j\right), x^* - \tilde x_j \rangle ) \notag\\
        & \hspace{1em}+\delta_2 \displaystyle \sum_{j = 0}^k A_j +  \delta_1 \displaystyle \sum_{j = 1}^{k} \alpha_j \|\tilde x_j - z_{j - 1} \|_2 \notag\\
        & \leqslant \frac{R^2}{2}  + A_k f(x^*) +  \delta_1 \displaystyle \sum_{j = 0}^{k} \alpha_j \|\tilde{x}_j - x^* \|_2  + \delta_2 \displaystyle \sum_{j = 0}^k A_j \notag\\ 
        & \hspace{1em} + \delta_1 \displaystyle \sum_{j = 1}^{k} \alpha_j \|\tilde x_j - z_{j - 1} \|_2.   \label{eq:stopping_rule_proof_1}
    \end{align}
    Whence, 
    \begin{align}
        \frac{1}{2} \|z_k - x^* \|_2^2 
        &\leqslant  \frac{R^2}{2} + A_k \Bigg( f(x^*) - f(x_k) +  \frac{\delta_1 \alpha_0}{A_k} \|\tilde{x}_0-x^*\|_2 \notag \\
        &  \hspace{1em} + \frac{\delta_1}{A_k} \displaystyle \sum_{j = 1}^{k} \alpha_j (2\|\tilde{x}_j - x^* \|_2 + \|z_{j-1} - x^* \|_2) +\frac{\delta_2}{A_k }\sum_{j = 0}^k A_j\Bigg). \label{eq:stopping_rule_proof_2}
    \end{align}
    Setting $k=0$, since $\|\tilde{x}_0-x^*\|_2 \leqslant R$ and, by the Theorem assumption, inequality \eqref{eq:stopping_rule_p} does not hold for $k\leqslant N-1$, we obtain
    \begin{align*}
    \frac{1}{2} \|z_0 - x^* \|_2^2 &\leqslant  \frac{R^2}{2} + A_0 \Bigg( f(x^*) - f(x_0) +  \frac{\delta_1 \alpha_0}{A_0} R + \delta_2 \Bigg) \\
    & \leqslant  \frac{R^2}{2} + A_0 \Bigg( - \delta_2 - 3R\delta_1 -\varepsilon +  \frac{\delta_1 \alpha_0}{A_0} R + \delta_2 \Bigg) \leqslant\frac{R^2}{2},
    \end{align*}
    where we also used that $\alpha_0=A_0$. Thus, we obtain that $\|z_0 - x^* \|_2\leqslant R$, and, since $x_0=z_0$, that $\|x_0 - x^* \|_2\leqslant R$. Hence, $\tilde{R}_0 \leqslant R$. Let us now assume that for some $k \leqslant N-1$, $\tilde{R}_{k-1} \leqslant R$ (see \eqref{equationR} for the definition of $\lbrace \tilde{R}_k \rbrace$). Then, by the definition of $  \tilde{x}_k  $ in Algorithm \ref{alg3} and convexity of the norm, we have that $\|\tilde{x}_k-x^*\|_2 \leqslant R$. Further, since $k \leqslant N-1$, we have that inequality \eqref{eq:stopping_rule_p} does not hold. Thus, from \eqref{eq:stopping_rule_proof_2}, we have:
    \begin{align*}
    \frac{1}{2} \|z_k - x^* \|_2^2& \leqslant  \frac{R^2}{2} + A_k \Bigg( f(x^*) - f(x_k) + \frac{\delta_1 \alpha_0 R}{A_k}  + \frac{\delta_1}{A_k} \displaystyle \sum_{j = 1}^{k} \alpha_j \cdot 3R +\frac{\delta_2}{A_k }\sum_{j = 0}^k A_j\Bigg) \\
    & \leqslant  \frac{R^2}{2} + A_k \Bigg( -  \frac{\delta_2}{A_k }\sum_{j = 0}^k A_j  - 3R\delta_1 -\varepsilon + 3R \delta_1 +\frac{\delta_2}{A_k }\sum_{j = 0}^k A_j \Bigg) \leqslant\frac{R^2}{2}.
    \end{align*}
    This implies that  $\|z_k - x^* \|_2\leqslant R$, and, by the definition of $x_k$ and the convexity of the norm, that $\|x_k - x^* \|_2\leqslant R$. Hence, $\tilde{R}_{k} \leqslant R$. In summary, we obtain by induction that for all $k \in \{0,\ldots,N-1\}$, $\tilde{R}_{k} \leqslant R$. This also implies that $\|\tilde{x}_N - x^* \|_2\leqslant R$.
    
    We now prove the second statement of the Theorem. Let us assume the opposite, i.e., $N > N_{\max}$. We use \eqref{eq:stopping_rule_proof_1} with $k=N-1$ and obtain, since $\tilde{R}_{N-1} \leqslant R$, that
    \begin{align*}
         f(x_{N-1}) -f(x^*)& \leqslant \frac{R^2}{2A_{N-1}} + 3R \delta_1  + \frac{\delta_2}{A_{N-1}} \displaystyle \sum_{j = 0}^{N-1} A_j  \\ &\leqslant \frac{2LR^2}{N^2} + 3R \delta_1  + \frac{\delta_2}{A_{N-1}} \displaystyle \sum_{j = 0}^{N-1} A_j \\
         &< \varepsilon + 3R \delta_1  + \frac{\delta_2}{A_{N-1}} \displaystyle \sum_{j = 0}^{N-1} A_j,
    \end{align*}
    where we used Claim \ref{cl3.3} and that $N>N_{\max}$. Thus, we see that after $N-1$ iterations, inequality \eqref{eq:stopping_rule_p} holds. This is a contradiction with the definition of $N$ as the first iteration number for which this inequality holds. This finishes the proof.
    \end{proof}
    
    Combining \eqref{eq:stopping_rule_p} with Claim \ref{cl3.4} and the fact that $N \leqslant N_{\max}$, we obtain that 
    \begin{equation*}
            f(x_N)-f(x^*) \leqslant \delta_2(N_{\max}+1) + 3R \delta_1 + \varepsilon.
    \end{equation*}
    Thus, if we redefine $\varepsilon \to \frac{\varepsilon}{3}$, and set $\delta_2 \leqslant \frac{\varepsilon}{3(N_{\max}+1)}$, $\delta_1 \leqslant \frac{\varepsilon}{9R}$, we guarantee that $f(x_N)-f(x^*) \leqslant \varepsilon$. 
    
    \begin{remark}
        \label{remark f star}
        In some situations we have at our disposal the value of $f(x^*)$ or its estimate. For example, when solving systems of linear equations by reformulating them as minimization problems:
        \begin{equation*}
            \begin{gathered}
                Ax = b, \\
                \min\limits_{x} \left\{f(x) = \|Ax - b \|_2^2 \right\},
            \end{gathered}
        \end{equation*}
        if a solution \pd{exists}, we have $f^* = f(x^*) = 0$.
        This allows us, based on \eqref{eq:stopping_rule_proof_1}, to change the inequality \eqref{eq:stopping_rule_p} to a more adaptive version, which can be checked online and which can be fulfilled much earlier than  \eqref{eq:stopping_rule_p}. Such counterpart of \eqref{eq:stopping_rule_p} reads as
        \begin{equation}
            \label{eq:stopping_rule_p_1}
            f(x_N)-f(x^*) \leqslant \frac{\delta_2}{A_N}\sum_{j=0}^{N}A_j + R \delta_1 + \delta_1 \displaystyle \sum_{j = 1}^{\pd{N}} \alpha_j \|\tilde x_j - z_{j - 1} \|_2 + \varepsilon.
        \end{equation}
        If this inequality is not fulfilled at iteration $k$, we have that $\tilde{R}_k \leqslant R$. If it is fulfilled at iteration $k$, we obtain that 
        \[
        f(x_k)-f(x^*) \leqslant \delta_2(k+1) + 3R \delta_1 + \varepsilon.
        \]
        Moreover, we also obtain that \eqref{eq:stopping_rule_p_1} holds after no more than $N_{\max}$ iterations. 
    \end{remark}

\section{Relative noise in the gradient}
In this section, we consider problem \eqref{eq:pr_st} in the relative noise setting (see \eqref{relative_inexact}), i.e., we assume that the inexact gradient $\tilde{\nabla} f(x)$ satisfies uniformly in $x \in Q$  the inequality
\begin{equation*}
 \|\widetilde{\nabla}f(x) - \nabla f(x) \| \leqslant \alpha \|\nabla f(x) \|_2.
\end{equation*}
As in the previous section, we assume that $f$ is $L_f$-smooth.
We also assume that $f$ is strongly convex with $\mu\ne 0$ and that $Q=\mathbb{R}^n$.
For this setting, we analyze a slightly different version of accelerated gradient method, adopted from \cite{stonyakin2021inexact}.

     \begin{algorithm}[H]
\caption{STM2 $(L, \mu_{\tau}, x_{start})$, $Q \subseteq \mathbb{R}^n$}
	\label{alg_from_adaptive}
\begin{algorithmic}[1]
\State 
\noindent {\bf Input:} Starting point $x_{start}$, number of steps $N$
\State {\bf Set} $y_0 = u_0 = x_0 = x_{start}$,
\State {\bf Set} $A_0 = \frac{1}{L}$, $\alpha_0 = A_0$.
\For {$k = 1 \dots N$}
        \State Find $\alpha_k$ from $(1 + \mu_{2} A_{k - 1})(A_{k - 1} + \alpha_k) = L \alpha_k^2$, \\
        \State or equivalently $\alpha_k = \frac{1 + \mu_{2} A_{k - 1}}{2L} + \sqrt{\frac{(1 + \mu_{2} A_{k - 1})^2}{4L^2} + \frac{A_{k - 1}(1 + \mu_{2} A_{k - 1})}{L}}$,
        \State $A_k = A_{k - 1} + \alpha_k,$
        \State $y_k = \frac{A_{k - 1} x_{k - 1} + \alpha_k u_{k - 1}}{A_k},$
        \\
        \State $\phi_k(x) = \alpha_k \langle \widetilde{\nabla}f(y_k), x - y_k \rangle + 
        \frac{1 + \mu_{2} A_{k- 1}}{2} \|u_{k - 1} - x \|_2^2 + \frac{\mu_{2} \alpha_k}{2} \|y_k - x \|_2^2,  $
        \\
        \State $u_k = \argmin_{u \in Q} \phi_k(u),$ 
        \State $x_k = \frac{A_{k - 1} x_{k - 1} + \alpha_k u_k}{A_k}.$
\EndFor
\State 
\noindent {\bf Output:} $x_N$.
\end{algorithmic}
\end{algorithm}
Since $Q = \mathbb{R}^n$, the main step of the algorithm can be simplified to
\begin{equation*}
    u_{k + 1} = \frac{1 + \mu_2 A_k}{1 + \mu_2 A_{k + 1}} u_k  + \frac{\mu_2 \alpha_{k + 1}}{1 + \mu_2 A_{k + 1}} y_{k + 1} - \frac{\alpha_{k + 1}}{1 + \mu_2 A_{k + 1}} \widetilde{\nabla} f(y_{k + 1}).
\end{equation*}

Combining Definition 3.3 of \cite{stonyakin2021inexact} with Claims~\ref{cl2.1}, \ref{cl2.33} and particular choice $V[y](x) = \frac{1}{2}\|x-y\|_2^2$, we have that $\delta$ in Definition 3.3 of \cite{stonyakin2021inexact} can be set to $\delta = \frac{3}{2}\frac{\delta^2}{\mu_2} \geqslant \frac{\delta^2}{2L_f} + \frac{\delta^2}{\mu} = \delta_2+\delta_3$, where we used that $\mu \leqslant L_f$ and that $\mu_2=\mu/2$. Further, $L$ in Definition 3.3 of \cite{stonyakin2021inexact} can be set to $L$ in our paper, and $\mu$ in Definition 3.3 of \cite{stonyakin2021inexact} can be set to $\mu_2=\mu/2$ in our paper. 
Algorithm \ref{alg_from_adaptive} is a particular case of Algorithm 2 in \cite{stonyakin2021inexact}. Since in this section, we are in the setting of relative inexactness \eqref{relative_inexact}, in each iteration of this algorithm we have a different error $\delta_k = \alpha\|\nabla f(y_k)\|_2$, which gives the following expression for $\delta_k$ in Algorithm 2 in \cite{stonyakin2021inexact}: $\delta_k = \frac{3\alpha^2\|\nabla f(y_k)\|_2^2}{2\mu_2}$.

Applying Theorem 3.4 of \cite{stonyakin2021inexact}, we obtain the following convergence rate for all $N\geqslant 0$:
\begin{equation*}
    \begin{gathered}
        f(x_N) - f(x^*) \leqslant \frac{R^2}{A_N} +\sum_{k = 1}^{N} 
        \frac{3\alpha^2 A_k \| \nabla{f}(y_k) \|_2^2}{\mu_2 A_N}:=\frac{\kappa}{A_N}, \\
        \|u_N - x^* \|_2^2 \leqslant \frac{1}{1 +A_N \mu_2} \left[ R^2 + \sum_{k = 1}^{N} 
        \frac{3\alpha^2 A_k \| \nabla{f}(y_k) \|_2^2}{\mu_2} \right]:=\frac{\kappa}{1 +A_N \mu_2}.
    \end{gathered}
\end{equation*}
Since we assumed that $Q=\mathbb{R}^n$, we have that $\nabla{f}(x^*) = 0$ and that, for all $x \in Q$, $ f(x) - f(x^*) \leqslant \frac{L}{4} \|x - x^* \|_2^2$, where we used \eqref{eq2.3:ref} and our definition $L=2L_f$.
Then, using convergence rate for $\lbrace u_k \rbrace$, we obtain
\begin{equation*}
    f(u_k) - f(x^*) \leqslant \frac{L}{4} \|u_k - x^* \|_2^2 \leqslant \frac{L\kappa}{4(1+ A_k \mu_2)}.
\end{equation*}
%
Using the convexity of $f$ and the definition of the sequence $\lbrace y_k \rbrace$ we get:
\begin{align*}
    f(y_{N + 1}) - f(x^*) & \leqslant \frac{\alpha_{N + 1}}{A_{N + 1}} \left[ f(u_N) - f(x^*) \right] + 
    \frac{A_N}{A_{N + 1}} \left[ f(x_N) - f(x^*) \right] \\
    & \leqslant \frac{\alpha_{N + 1}}{A_{N + 1}} \cdot \frac{L\kappa}{4(1+ A_k \mu_2)} + \frac{A_N}{A_{N + 1}} \cdot \frac{\kappa}{A_N}. 
\end{align*}
Our next goal is to estimate $\frac{\alpha_{N + 1}}{A_{N + 1}} \cdot \frac{L}{4(1+ A_k \mu_2)} $ from above.
Using the inequalities $    \frac{A_k}{1 + \mu_2 A_k} \leqslant \frac{1}{\mu_2}$ and $\sqrt{x + y} \leqslant \sqrt{x} + \sqrt{y}$, and the definition of the sequence $\lbrace \alpha_k \rbrace$:
\begin{equation*}
    \begin{gathered}
        \alpha_k = \frac{1 + \mu_2 A_{k - 1}}{2L} + \sqrt{\frac{(1 + \mu_2 A_{k - 1})^2}{4L^2} + \frac{A_{k - 1}(1 + \mu_2 A_{k - 1})}{L}}, \\
    \end{gathered}
\end{equation*}
we have
\begin{multline*}
    \frac{\alpha_{N + 1}}{A_{N + 1}} \cdot \frac{L}{4(1+ A_k \mu_2)} \\ =\frac{L}{4A_{N + 1}} \frac{1}{1 + \mu_2 A_{N}} \left( \frac{1 + \mu_2 A_{N}}{2L} + \sqrt{\frac{(1 + \mu_2 A_{N})^2}{4L^2} + \frac{A_{N}(1 + \mu_2 A_{N})}{L}} \right) \\ \leqslant
    \frac{L}{4}\frac{1}{A_{N + 1}} \left(\frac{1}{2L} + 
    \frac{1}{2L} + \sqrt{\frac{1}{L\mu_2}} \right) \leqslant
    \frac{1}{4A_{N + 1}} \sqrt{\frac{L}{\mu_2}}.
\end{multline*}
This gives us the following estimate
\begin{multline*}
    f(y_{N+1}) - f^* \leqslant \frac{\kappa}{4A_{N + 1}} \sqrt{\frac{L}{\mu_2}} + \frac{\kappa}{A_{N+1}}
    \leqslant
    \frac{1}{4} \sqrt{\frac{L}{\mu_2}} \left( \frac{5R^2}{A_{N + 1}} + \sum_{k = 1}^{N} \frac{15 \alpha^2 A_k \| \nabla{f}(y_k) \|_2^2}{ \mu_2 A_{N + 1}} \right),    
\end{multline*}
where we used that $L/\mu_2 \geqslant 1$ and the definition of $\kappa$.

Since $f$ is $L_f$-smooth, $L=2L_f$ and $\nabla f(x^*)=0$, we obtain for any $x\in Q$ that
\begin{equation*}
    \| \nabla{f}(x) \|_2^2 \leqslant L \left( f(x) - f(x^*) \right).
\end{equation*}
Whence, using the previous bound,
\begin{equation*}
    f(y_{N + 1}) - f(x^*) \leqslant \frac{1}{4} \sqrt{\frac{L}{\mu_2}} \left( \frac{5R^2}{A_{N + 1}} + \sum_{k = 1}^{N} \frac{15 \alpha^2 A_k L(f(y_k) - f^*)}{\mu_2 A_{N + 1}} \right).
\end{equation*}
Introducing the following notations $\lambda = \frac{5R^2}{4} \sqrt{\frac{L}{\mu_2}}$, $\theta = \frac{15\alpha^2}{4} \sqrt{\frac{L^3}{\mu_2^3}}$, $\Delta_k = f(y_k) - f(x^*)$, we obtain the following recurrence
%
\begin{equation*}
    \Delta_{N} \leqslant \frac{\lambda}{A_N} + \theta \sum_{k = 0}^{N - 1}\frac{A_k}{A_N} \Delta_k,
\end{equation*}
where we add the term corresponding to $k=0$ to the sum to simplify the proof that will follow. Analyzing this recurrence, we obtain.
\begin{claim} \label{claim_ind}
For all $k\geqslant 1$ it holds that 
\begin{equation*}
    \Delta_k \leqslant \frac{\left(1 + \theta\right)^{k - 1}}{A_k}\lambda + \theta\frac{A_0 \left(1 + \theta\right)^{k - 1}}{A_k}\Delta_0.
\end{equation*}
\end{claim}

\begin{proof}
    The induction basis $k = 1$ is obvious.
    Induction step:
    \begin{align*}
        \Delta_k& \leqslant \frac{\lambda}{A_k} + \theta \displaystyle\sum_{j = 0}^{k - 1} \frac{A_j}{A_k}\Delta_j \\
        & \leqslant \frac{\lambda}{A_k} + \theta \sum_{j = 1}^{k - 1} \frac{A_j}{A_k} \Delta_j + 
        \theta \frac{A_0}{A_k} \Delta_0 \\
        &\leqslant \frac{\lambda}{A_k} + \theta \displaystyle\sum_{j = 1}^{k - 1} \left( \frac{A_j}{A_k} \frac{(1 + \theta)^{j - 1}}{A_k} \lambda + \theta \frac{A_0(1 + \theta)^{j - 1}}{A_k} \Delta_0 \right)  + \frac{A_0}{A_k}\Delta_0 \\
        &\leqslant \frac{\lambda}{A_k} + \theta \displaystyle\sum_{j = 0}^{k - 2}\left(\frac{\lambda (1 + \theta)^j}{A_k} + \theta \frac{A_0 (1 + \theta)^j}{A_k} \Delta_0 \right) + \frac{A_0}{A_k}\Delta_0 \\
        & = \frac{1}{A_k} \left( \lambda + \lambda \left[(1 + \theta)^{k - 1} - 1 \right] + 
        \theta A_0 \Delta_0 \left[ (1 + \theta)^{k - 1} - 1\right] + A_0 \Delta_0 \right) \\
        & = \frac{\left(1 + \theta\right)^{k - 1}}{A_k}\lambda + \theta\frac{A_0 \left(1 + \theta\right)^{k - 1}}{A_k}\Delta_0.
    \end{align*}
\end{proof}
This gives us the following result
\begin{equation*}
    f(y_k) - f(x^*) \leqslant \frac{\lambda(1 + \theta)^k}{A_k} + \theta \frac{A_0(1 + \theta)^k}{A_k} \left(f(y_0) - f(x^*) \right).
\end{equation*}
By the definition of $\theta = \frac{15\alpha^2}{4} \sqrt{\frac{L^3}{\mu_2^3}}$, we obtain, that, if we choose $\alpha$ as
\begin{equation} \label{alpha_bound}
    \begin{gathered}
        \alpha \leqslant \frac{1}{7} \frac{\mu_2}{L} = \Theta \left( \frac{\mu}{L} \right),
    \end{gathered}
\end{equation}
then
\begin{equation*}
    \begin{gathered}
    1 + \sqrt{\frac{\mu_2}{L}} \left(\frac{1}{2} + \frac{15}{196} + \frac{15}{392} \right) \leqslant 1
        + \sqrt{\frac{\mu_2}{L}} \sqrt{\frac{1}{2}}
        \\
        \Leftrightarrow
        1 + \theta + \frac{1}{2}\sqrt{\frac{\mu_2}{L}} +
        \frac{1}{2}\sqrt{\frac{\mu_2}{L}} \theta \leqslant
        1 + \sqrt{\frac{\mu_2}{2L}}
        \\
        \Leftrightarrow
        \frac{1 +\theta}{1 + \sqrt{\frac{\mu_2}{2L}}} \leqslant \frac{1}{1 + \frac{1}{2}\sqrt{\frac{\mu_2}{L}}}.
    \end{gathered}
\end{equation*}
Combining this with Claim \ref{cl3.1} and Corollary~\ref{lower bound lambda}, we obtain that
\begin{equation} \label{eq_bound_fin}
 \frac{(1 +\theta)^k}{A_k} \leqslant  \left(\frac{1 +\theta}{1 + \sqrt{\frac{\mu_2}{2L}}} \right)^k \frac{1}{A_0} \leqslant \left(\frac{1}{1 + \frac{1}{2}\sqrt{\frac{\mu_2}{L}}} \right)^k \frac{1}{A_0} \leqslant L \exp{\left({\displaystyle-\frac{k}{4}\sqrt{\frac{\mu_2}{L}}} \right)}.
\end{equation}
As a result, we get the following theorem.
\begin{theorem} \label{th_alph}
    Assume that the objective $f$ is $L_f$-smooth and strongly convex with $\mu\ne 0$, that the inexactness  in the gradient is described by \eqref{relative_inexact}, and that $Q=\RR^n$. Also assume that $\alpha$ is chosen according to \eqref{alpha_bound}. Then, for all $k\geqslant 1$, the sequence $\lbrace y_k \rbrace$ generated by Algorithm \ref{alg_from_adaptive} satisfies
    \begin{equation*}
        f(y_k) - f(x^*) \leqslant \left(\frac{5LR^2}{4} + \frac{15}{196} \sqrt{\frac{2L}{\mu}}\left[f(y_0) - f(x^*) \right] \right)\exp{\left({\displaystyle-\frac{k}{4}\sqrt{\frac{\mu}{2L}}} \right)}.
    \end{equation*}
\end{theorem}
%

\section{\pd{Extensions}}
\label{section_extinsons}
In this section, we extend the analysis of Algorithm \ref{alg3} with absolute noise to two settings. The first extension is an extension to stochastic optimization setting where the error in the gradient has stochastic nature. The second one is the extension to structured nonsmooth setting of composite minimization, where the objective is given as a sum of smooth part with inexact gradient and a simple convex function. In both cases, the analysis mainly follows the lines of Section \ref{section 4}. Thus, we underline the differences and skip in the proofs some steps that are similar to the analysis in that section.

\subsection{Random additive noise in the gradient}
\label{stoh_optim_section}
In this subsection, we extend the  analysis of Algorithm~\ref{alg3} for the setting of random absolute noise in the gradient. 
We assume that an algorithm can use the stochastic gradient $\widetilde{\nabla} f(x, \xi)$, 
which is assumed to have bounded variance for all, possibly random, $x\in Q$:
\begin{equation}
    \label{stohastic_error_rand}
    \mathbb{E}_{\xi} \left[ \|\widetilde{\nabla}f(x, \xi) - \nabla{f}(x) \|_2^2 \; \Big| \; x \right] \leqslant \delta^2.
\end{equation}
Similarly to Section \ref{section 4}, we assume $L_f$-smoothness and $\mu$-strong convexity of $f$, i.e., that~\eqref{eq2.3:ref},\eqref{eq2.4:ref} hold.
As before, we set $L=2L_f$ and 
\begin{equation*}
    \begin{gathered}
        \delta_1 = \delta,   \quad
        \delta_2 = \frac{\delta^2}{L},  \quad
        \delta_3 = \frac{\delta^2}{\mu},
    \end{gathered}
\end{equation*}
where the latter quantity is defined whenever $\mu >0$.

One of the main motivations for such stochastic problems is machine learning. For example, Empirical Risk Minimization problem with the finite-sum structure of the objective
\begin{equation*}
    f(x) = \frac{1}{M}\sum_{i = 1}^M f_i(x)
\end{equation*}
can be considered as a stochastic optimization problem with stochastic gradient
\begin{equation*}
    \begin{gathered}
        \widetilde{\nabla} f(x, \xi) = \frac{1}{m}\sum_{i \in \xi} \nabla f_i(x); \; \xi \subset \lbrace 1 \dots M \rbrace, \; |\xi| = m, \; m < M,
    \end{gathered}
\end{equation*}
where $\xi$ is a random subset of $\lbrace 1 \dots M \rbrace$.
It should be noted that the error $\delta^2$ can be reduced by the use of mini-batches. Namely increasing the size of $\xi$ from 1 to $m$ decreases the variance from $\delta^2$ to $\frac{\delta^2}{m}$.

The first step of the analysis is to obtain the counterparts of Claims~\ref{cl2.1}, \ref{cl2.2}, and \ref{cl2.33} in the stochastic setting. 
\begin{claim}
\label{claim_stoh_L}
Assume that $x, y$ are random vectors. Then, 
\begin{equation*}
    \mathbb{E} \left[ f(y) \right] \leqslant \mathbb{E} \left( f(x) + \langle \widetilde{\nabla}f(x, \xi), y - x \rangle + \frac{L}{2} \|x - y \|_2^2 \right) + \delta_2.
\end{equation*}
\end{claim}

\begin{proof}
Using the $L_f$-smoothness, we obtain
    \begin{equation*}
        \begin{split}
            \mathbb{E} \left[ f(y) \; \Big| \; x \right] &
            \leqslant \mathbb{E} \left[f(x) + \langle \nabla f(x), y - x \rangle + \frac{L_f}{2} \|y - x \|_2^2  \; \Big| \; x \right] \\
            & \stackrel{\eqref{stohastic_error_rand}}{\leqslant} \mathbb{E} \left[f(x) + \langle \widetilde{\nabla} f(x, \xi), y - x \rangle + \frac{L}{2} \|y - x \|_2^2 \; \Big| \; x \right] + \delta_2,
        \end{split}
    \end{equation*}
where $L=2L_f$.
Taking the full expectation of both sides, we get the required.
\end{proof}
Using the same steps as in the proof of Claim~\ref{claim_stoh_L}, we get the counterparts of Claims~\ref{cl2.2}, \ref{cl2.33}.
\begin{claim} 
\label{claim_mu_1_stoh}
Assume that $x, y$ are random vectors. Then,
\begin{equation*}
    \mathbb{E} \left[f(x) + \langle \widetilde{\nabla}f(x, \xi), y - x \rangle + \frac{\mu}{2} \|x - y \|_2^2 - \delta_1\|x - y \|_2 \; \Big| \; x \right] \leqslant \mathbb{E} \left[f(y) \; \Big| \; x \right].
\end{equation*}
\end{claim}
\begin{claim} \label{claim_mu_2_stoh}
Assume that $x, y$ are random vectors and that \(\mu \not= 0 \).  Then, 
\begin{equation*}
    \mathbb{E} \left[f(x)+ \langle \widetilde{\nabla} f(x, \xi), y - x \rangle + \frac{\mu}{4} \|y - x \|_2^2 \; \Big| \; x \right] - \delta_3 \leqslant \mathbb{E} \left[f(y) \; \Big| \; x \right].
\end{equation*}
\end{claim}
The following sequence is the counterpart of the sequence $\lbrace \tilde{R}_k\rbrace$:
\begin{equation}\label{eq:tildeBN_def}
    \tilde{B}_k = \max\limits_{0 \leqslant j \leqslant k} \lbrace \mathbb{E}\|z_j - x^*\|_2, \mathbb{E}\|x_j - x^*\|_2, \mathbb{E}\|\tilde x_j - x^*\|_2 \rbrace.
\end{equation}
Using the above, we obtain the following counterparts of \pd{Propositions}~\ref{thr4.1} and \ref{th4.2} under the assumptions of this subsection.
\begin{proposition} \label{th_stoh_1}
    The sequences generated by Algorithm~\ref{alg3} satisfy for all $k \geqslant 0$ the
    inequality:
    \begin{equation*}
        A_k \mathbb{E} \left[f(x_k)\right] \leqslant \mathbb{E} \left[ \psi_k(z_k) \right] + \delta_2 \displaystyle\sum_{j=0}^{k} A_j + \delta_1 \displaystyle\sum_{j=1}^{k} \alpha_j \mathbb{E}\left[\|\tilde{x}_j-z_{j-1}\|_2 \right].
    \end{equation*}
\end{proposition}
\begin{proposition} \label{th_stoh_2}
    If $ \mu \not= 0 $, the sequences generated by Algorithm~\ref{alg3} satisfy for all $k \geqslant 0$ the
    inequality:
    \begin{equation*}
        A_k \mathbb{E} \left[f(x_k)\right] \leqslant \mathbb{E} \left[\psi_k(z_k) \right] +\delta_2 \displaystyle\sum_{j=0}^{k} A_j + \delta_3 \displaystyle\sum_{j=0}^{k - 1} A_j. 
    \end{equation*}
\end{proposition}
The proofs of these propositions repeat the same induction steps as in the proofs of Propositions~\ref{thr4.1}, \ref{th4.2}, but using the new Claims~\ref{claim_stoh_L}, \ref{claim_mu_1_stoh}, and \ref{claim_mu_2_stoh}. Using the last two propositions, we finally  obtain the following counterpart of convergence Theorem~\ref{th4conv} for the stochastic setting.
\begin{theorem}[Convergence rate of stochastic STM] \label{th_conv_stoh}
$\\$
Let $ \|\tilde{x}_0 - x^* \|_2 \leqslant R$ for some $R$, function $f$ be $L_f$-smooth and strongly convex with parameter $\mu \geqslant 0$. Let the stochastic gradient $\widetilde{\nabla}f(x, \xi)$ satisfy
\begin{equation}
    \begin{gathered}
        \mathbb{E}_{\xi} \left[\|\widetilde{\nabla}f(x, \xi) - \nabla{f}(x) \|_2^2 \Big| \; x \right]  \leqslant \delta^2.
    \end{gathered}
\end{equation}
Then, if $ \mu \not= 0 $, the sequence $\lbrace x_N \rbrace$ generated by Algorithm $\ref{alg3}$ satisfy for all $N \geqslant 0$ the inequalities
\begin{align*}
        \mathbb{E} f(x_N) - f(x^*) &\leqslant L R^2 \exp\left({{-\frac{1}{2}\sqrt{\frac{\mu_1}{L}}N}}\right) + \left(1 + \sqrt{{\frac{L}{\mu_1}}}\right)\delta_2 + 3 \tilde B_N \delta_1, \\ 
        \mathbb{E} f(x_N) - f(x^*) & \leqslant L R^2 \exp{\left(-\frac{1}{2}\sqrt{\frac{\mu_2}{L}}N \right)} + \left(1 + \sqrt{{\frac{L}{\mu_2}}} \right)\delta_2 + \left(1 + \sqrt{\frac{L}{\mu_2}} \right)\delta_3. \\
\end{align*}
If $ \mu = 0 $, the sequences $\lbrace x_k \rbrace$, $\lbrace \tilde{x}_k \rbrace$, $\lbrace z_k \rbrace$ generated by Algorithm $\ref{alg3}$ satisfy for all $N \geqslant 0$ the inequality   
\begin{equation*}
        \mathbb{E} f(x_N) - f(x^*) \leqslant \frac{4LR^2}{N^2} + 3 \tilde B_N \delta_1 + N \delta_2.
\end{equation*}
Here the sequence $\lbrace \tilde{B}_k \rbrace$ is defined in \eqref{eq:tildeBN_def}.
\end{theorem}
As we see, Algorithm~\ref{alg3} has the same convergence rate in the stochastic setting as in the deterministic setting. The proof of the above theorem repeats the same steps as the proof of Theorem~\ref{th4conv}. Thus, we omit the proof.
\begin{remark}
    Usually, in the context of stochastic optimization, the analysis of algorithms relies also on the assumption of unbiased stochastic gradient:
    \begin{equation*}
        \mathbb{E} \left[\widetilde{\nabla} f(x, \xi) \; \Big| \; x \right] = \nabla{f}(x).
    \end{equation*}
    Our analysis does not require this assumption.
\end{remark}


\subsection{Nonsmooth objective}
\label{Nonsmooth_structure}
In this subsection, we consider the problem of structured nonsmooth optimization, usually referred to as composite minimization, 
    \begin{equation}
    \label{eq:pr_st_composite}
       \min_{x \in Q} \left\{ f(x) = \mathcal{L}(x) + r(x)\right\}.
    \end{equation}
    We assume that the function $\mathcal{L}$ is $L_f$-smooth and $\mu$-strongly-convex (see $\eqref{eq2.2:ref}, \eqref{eq2.4:ref}$), the function $r(x)$ is convex and relatively simple. We further assume that inexact gradient $\tilde{\nabla} \mathcal{L}(x)$ with absolute noise (cf. \eqref{inexact}) is available for $\mathcal{L}$.
    
    This setting is motivated, in particular, by machine learning problems, for example, logistic regression loss minimization problem with the $l_1$ regularization and dataset $\lbrace(X_i, y_i)\rbrace_{i = 1}^{K}$, where \pd{$y_i \in \lbrace 0, 1 \rbrace$ for $i \in \lbrace 1, K \rbrace$}. For this problem, we have
    \begin{equation*}
        \begin{gathered}
            \mathcal{L}(x) = \sum_{i = 1}^{K} y_i \ln{p_i(x)} + (1 - y_i) \ln{(1 - p_i(x))}, \\
            p_i(x) = \sigma(\langle x, X_i \rangle), 
            \sigma(z) = \frac{1}{1 + \exp{(-z)}}, \\
            r(x) = \lambda_1 \|x \|_1.
        \end{gathered}
    \end{equation*}
    
    In the setting of composite minimization, Algorithm $\ref{alg3}$ requires only one change in the definition of the function sequence $\lbrace \psi_k(\cdot) \rbrace$ as follows:
    \begin{equation}
    \label{eq:comp_psi_k_def}
        \begin{gathered}
            \psi_0(x) = \frac{1}{2}\|x - \tilde{x}_0\|_2^2 + \alpha_{0} \left( \mathcal{L}(\tilde{x}_0) + \langle \tilde{\nabla} \mathcal{L}(\tilde{x}_0), x - \tilde{x}_0 \rangle + \frac{\mu_{\tau}}{2} \|x - \tilde{x}_0 \|_2^2 + r(x) \right), \\
            \psi_k(x) = \psi_{k - 1}(x) + \alpha_{k} \left( \mathcal{L}(\tilde{x}_k) + \langle \tilde{\nabla} \mathcal{L}(\tilde{x}_k), x - \tilde{x}_k \rangle + \frac{\mu_{\tau}}{2} \|x - \tilde{x}_k \|_2^2 + r(x) \right).
        \end{gathered}
    \end{equation}
    For such modified algorithm, in the concept of absolute noise~\eqref{inexact}, the convergence result remains the same. However, some intermediate statements, such as Lemma~\ref{LEM}, require a different analysis.
    Therefore, we make a different analysis to obtain an estimate in the spirit of Proposition~\ref{thr4.1}.
    \begin{lemma}[auxiliary statement for $\psi_k$'s]
    \label{lemma for nonsmth}
    Under the assumptions of this subsection, for the modified sequence $\lbrace \psi_k(\cdot) \rbrace$, we have
    \begin{multline*}
        \psi_{k + 1}(z_{k + 1})
        \geqslant \psi_k(z_k) + \frac{1 + \mu_{\tau}A_k}{2}\|z_k - z_{k + 1}\|_2^2
        \\
        + \alpha_{k} \left( \mathcal{L}(\tilde{x}_k) + \langle \tilde{\nabla} \mathcal{L}(\tilde{x}_k), x - \tilde{x}_k \rangle + \frac{\mu_{\tau}}{2} \|x - \tilde{x}_k \|_2^2 + r(x) \right).
    \end{multline*}
    \end{lemma}
    \begin{proof}
        The function $\psi_k$ defined in \eqref{eq:comp_psi_k_def} is $\frac{1 + \mu_{\tau}A_k}{2}$-strongly-convex. Thus, since $z_k$ is its minimizer, we have 
        \begin{equation*}
            \psi_k(z_{k + 1}) \geqslant \psi_k(z_k) + \frac{1 + \mu_{\tau}A_k}{2}\|z_k - z_{k + 1}\|_2^2.
        \end{equation*}
        Using the recurrent definition of $\psi_{k + 1}$ we obtain the required by induction.
    \end{proof}
    Using Lemma~\ref{lemma for nonsmth} instead of Lemma~\ref{LEM} and convexity of the function $r(x)$, we can obtain a result similar to \pd{Proposition}~\ref{thr4.1}.
    \begin{proposition}
    \label{prop_constr_ind}
    The sequences $\lbrace x_k \rbrace$, $\lbrace \tilde{x}_k \rbrace$, $\lbrace z_k \rbrace$ generated by Algorithm $\ref{alg3}$ modified for structured nonsmooth optimization satisfy for all $k \geqslant 0$ the inequality
    \begin{equation*}
        A_k f(x_k) \leqslant \psi_k(z_k) + \delta_2 \displaystyle\sum_{j=0}^{k} A_j + \delta_1 \displaystyle\sum_{j=1}^{k} \alpha_j \|\tilde{x}_j-z_{j-1}\|_2.  
    \end{equation*}
    \end{proposition}
    \begin{proof}
        The induction basis $k = 0$ is obvious and repeats the proof of Proposition~\ref{thr4.1} since
        \begin{equation*}
            f(x_0) \leqslant L \psi_0(z_0) + \delta_2.
        \end{equation*}
        Let us consider iteration $k >0$. Since $r(x)$ is convex, by the definition of $ x_k$, we get:
        \begin{equation*}
            A_k r(x_k) \leqslant A_{k - 1} r(x_{k - 1}) + \alpha_{k} r(z_k).
        \end{equation*}
        By this inequality, Claim \ref{cl2.1} applied to $\mathcal{L}(x)$, the definition of the sequences $\lbrace  x_k \rbrace$, $\lbrace \tilde{x}_k \rbrace$,  we have
        \begin{align*}
            A_k f(x_k) &\leqslant 
            A_k \left( \mathcal{L}(\tilde x_{k}) + \langle \tilde \nabla \mathcal{L}(\tilde x_{k}), x_k - \tilde x_k \rangle + \frac{L}{2} \|x_k - \tilde x_k \|_2^2 + \delta_2 + r(x_k) \right) \\
            & \leqslant  A_{k - 1} \left( \mathcal{L}(\tilde x_k) + \langle \tilde \nabla \mathcal{L}(\tilde x_k), x_{k - 1} - \tilde x_k \rangle + r(x_{k - 1})\right) \\
            & \hspace{1em}+ \alpha_k \left( \mathcal{L}(\tilde x_k) + \langle \tilde \nabla \mathcal{L}(\tilde x_k), z_k - \tilde x_k \rangle + r(z_k) \right) + A_k \delta_2 + \frac{L \alpha_k^2}{A_k}\|z_k - z_{k - 1} \|_2^2\\
            & \leqslant A_{k - 1}f(x_{k - 1}) + \alpha_k \left( \mathcal{L}(\tilde x_k) + \langle \tilde \nabla \mathcal{L}(\tilde x_k), z_k - \tilde x_k \rangle + r(z_k) \right) \\
           & \hspace{1em}+  \frac{1 + \mu_{\tau} A_{k - 1}}{2}\|z_k - z_{k - 1} \|_2^2 + A_k \delta_2 + A_{k-1} \delta_1 \|x_{k - 1} - \tilde x_k \|_2,
        \end{align*}
        where in the last inequality we used the equation in Step \ref{step:alpha_k} of Algorithm $\ref{alg3}$ and Claim~\ref{cl2.2} applied to $\mathcal{L}$.
        By the induction hypothesis and since $\frac{\mu_{\tau}}{2}\|z_k-\tilde{x}_k\|_2^2 \geqslant 0$, we further obtain
        \begin{align*}
             A_k f(x_k) - A_{k-1} \delta_1 \|x_{k - 1} - \tilde x_k \|_2 
             &\leqslant \psi_{k - 1}(z_{k - 1}) + \delta_2 \displaystyle\sum_{j=0}^{k - 1} A_j + \delta_1 \displaystyle\sum_{j=1}^{k-1} \alpha_j \|\tilde{x}_j-z_{j-1}\|_2 \\ 
             &+ \alpha_k \left( \mathcal{L}(\tilde x_k) + \langle \tilde \nabla \mathcal{L}(\tilde x_k), z_k - \tilde x_k \rangle +\frac{\mu_{\tau}}{2}\|z_k-\tilde{x}_k\|_2^2 + r(z_k) \right) \\
             &+ \frac{1 + \mu_{\tau} A_{k - 1}}{2} \|z_k - z_{k - 1} \|_2^2  + A_k \delta_2.
        \end{align*}
        Using Lemma~$\ref{lemma for nonsmth}$, we can finish the proof in a similar way as in the proof of Proposition~\ref{thr4.1}.
    \end{proof}

    We finally obtain the following counterpart of Theorem~\ref{th4conv} for composite minimization problems.
    \begin{theorem}
    Let the modified Algorithm \ref{alg3} be applied to composite problem \eqref{eq:pr_st_composite}, 
    where the function  $\mathcal{L}(x)$ is $L_f$-smooth and $\mu$-strongly-convex and the function $r(x)$ is convex. If $ \mu \not= 0 $, the sequences $\lbrace x_k \rbrace$, $\lbrace \tilde{x}_k \rbrace$, $\lbrace z_k \rbrace$ generated by the modified Algorithm $\ref{alg3}$ satisfy for all $N \geqslant 0$ the inequalities
        \begin{equation*}
            \begin{gathered}
                f(x_N) - f(x^*) \leqslant L R^2 \exp\left({{-\frac{1}{2}\sqrt{\frac{\mu_1}{L}}N}}\right) + \left(1 + \sqrt{{\frac{L}{\mu_1}}}\right)\delta_2 + 3 \tilde R_N \delta_1, \\ 
                f(x_N) - f(x^*) \leqslant L R^2 \exp{\left(-\frac{1}{2}\sqrt{\frac{\mu_2}{L}}N \right)} + \left(1 + \sqrt{{\frac{L}{\mu_2}}} \right)\delta_2 + \left(1 + \sqrt{\frac{L}{\mu_2}} \right)\delta_3. \\
            \end{gathered}
        \end{equation*}
        If $ \mu = 0 $, the sequences $\lbrace x_k \rbrace$, $\lbrace \tilde{x}_k \rbrace$, $\lbrace z_k \rbrace$ generated by the modified Algorithm $\ref{alg3}$ satisfy for all $\pd{N} \geqslant 0$ the inequality   
        \begin{equation*}
                f(x_N) - f(x^*) \leqslant \frac{4LR^2}{N^2} + 3 \tilde R_N \delta_1 + N \delta_2,
        \end{equation*}
        where the sequence $\lbrace \tilde R_k \rbrace$ is defined in \eqref{equationR}.
    \end{theorem}
    As we see, for composite problems, modulo a small modification of the algorithm, the main result is the same as in the smooth case.


\section{Conclusions and observations}
In this section, we give a number of remarks in order to discuss the obtained results. \pd{In particular, the convergence rate results obtained so far explicitly include the oracle inexactness, and we can look at these results from a little bit different angle of controlling the inexactness. In particular, if the oracle error can be controlled, we can estimate how small should be the oracle error if our goal is to obtain an $\varepsilon$-solution to the problem. Such bound also give an estimate for the largest tolerable error not preventing the algorithm from obtaining an $\varepsilon$-solution.}

\pd{
\begin{remark}
    \label{remake future work}    
    In Sections~\ref{stoh_optim_section}, \ref{Nonsmooth_structure}, we considered the extensions of Algorithm \ref{alg3} with absolute noise to the settings of stochastic optimization and structured nonsmooth optimization.
    We strongly believe that it is possible to combine these two extensions into one since the analysis in both cases follows the same lines as the analysis in Section \ref{section 4}.
    We believe that the same can be also done with the analysis of Algorithm \ref{alg_from_adaptive} under the relative noise in the gradient (see stochastic version of this condition in \cite{vaswani2019fast}). We leave these developments for the future work.
\end{remark}
}
\begin{remark}
    \label{remark 5.1}
    The results of Theorem~\ref{th4conv} and Proposition \(\ref{th4.3} \) are obtained for possibly unbounded feasible set $Q$. If this set is compact, we can set $R = \mathrm{diam}(Q)$, i.e., the diameter of the set $Q$. This simplifies the results and derivations since, in this case, by the construction of Algorithm \ref{alg3},  \(\tilde R_k \leqslant R \) for all $k\geqslant 0$.
\end{remark}

%

\begin{remark}
    \label{remark 5.3}
    When considering the absolute noise, in Section \ref{sec2}, we had two possibilities for dealing with ``inexact strong convexity'': according to Claim  \ref{cl2.2} when $\mu \geqslant 0$ and according to Claim \ref{cl2.33} when $\mu > 0$. This resulted in two different bounds in Theorem~\ref{th4conv} in the setting when $\mu>0$. 
Recalling that
    \begin{equation*}
        \begin{gathered}
            \delta_1 = \delta, \quad
            \delta_2 = \frac{\delta^2}{L},  \quad
            \delta_3 = \frac{\delta^2}{\mu},
        \end{gathered}
    \end{equation*}
and comparing the two bounds in Theorem~\ref{th4conv}, we see that if
    \begin{equation*}
        \delta < \frac{3 \tilde{R}}{\frac{1 + \sqrt{\frac{L}{\mu}}}{\mu} + \frac{\sqrt{\frac{L}{\mu}}(\sqrt{2} - 1)}{L}},
    \end{equation*}
then the model corresponding to \(\tau = 2\), that is described in Claim \ref{cl2.33}, leads to a smaller term in the convergence rate bound due to the error accumulation than the model  corresponding to $\tau = 1$, that is described in Claim \ref{cl2.2}
\end{remark}

%

\pd{
The above results are valid for uncontrolled and unknown values of the error $\delta$ in the model of absolute noise. 
At the same time, in some cases, it may happen that the error $\delta$ can be controlled and made as small as one desires. For example, in the setting of Section \ref{S:inverse_problems}, the gradient can be approximated using finite-difference solution of primal and adjoint systems of equations, and $\delta$ can be decreased by decreasing the discretization step. In the setting of Section \ref{stoh_optim_section}, the error $\delta$ can be made smaller by the means of using mini-batches of stochastic gradients.
Thus, a natural question is how small should one choose the accuracy $\delta$ if the goal is to find an $\varepsilon$-approximate solution, i.e., guarantee $f(x_k) - f(x^*) \leqslant \varepsilon$? A similar question could be as follows: given a target accuracy $\varepsilon$, how large is the  error  $\delta$ that can be tolerated by an algorithm still guaranteeing the target accuracy $\varepsilon$? 
This, in particular, allows one to compare the robustness of different algorithms with respect to the noise.
In the following series of remarks, we address these questions by deriving the relations between $\delta$ and $\varepsilon$. 
}

\begin{remark}
    \label{remark5.4}
    Let us consider the ``inexact strong convexity'' model corresponding to $\tau = 2$, $\mu>0$, $\mu_2 = \frac{\mu}{2}$, and $\delta_2 = \frac{\delta^2}{2 L_f}$, $L = 2 L_f$, $\delta_3 = \frac{\delta^2}{\mu}$ (see Claims \ref{cl2.1}, \ref{cl2.33}). In this case, we can write explicit expressions for the dependence of the error $\delta$ and the iteration number $N$  on the target accuracy $\varepsilon$. Substituting the above values into the bound in Theorem~\ref{th4conv}, we obtain 
\begin{align*}
        f(x_N) - f(x^*) &\leqslant L R^2 \exp{\left(-\frac{1}{2}\sqrt{\frac{\mu}{2L}}N \right)} + \left(1 + \sqrt{{\frac{L}{\mu_2}}} \right)\delta_2 + \left(1 + \sqrt{\frac{L}{\mu_2}} \right)\delta_3 \\
        &=  2L_f R^2 \exp{\left(-\frac{1}{2}\sqrt{\frac{\mu}{4L_f}}N \right)} + \left(1+ \sqrt{\frac{4L_f}{\mu}} \right) \left( \frac{2L_f+\mu}{2\mu L_f} \right) \delta^2.\\
    \end{align*}
Thus, choosing
 \begin{equation*}
        \begin{gathered}
            \delta \leqslant \sqrt{\varepsilon} \sqrt{\frac{\mu L_f}{\mu + 2 L_f}} \left(  1+ \sqrt{\frac{4L_f}{\mu}} \right)^{-\frac12}
            \pd{
            \leqslant \left(\frac{L_f \varepsilon}{1 + 2 \frac{L_f}{\mu} + 2 \sqrt{\frac{L_f}{\mu}} + 4 \frac{L_f \sqrt{L_f}}{\mu \sqrt{\mu}}} \right)^{\frac{1}{2}}
            = O\left(\sqrt{\mu \varepsilon} \left(\frac{\mu}{L_f} \right)^{\frac{1}{4}} \right);}\\
            N \geqslant 2 \sqrt{\frac{4L_f}{\mu}} \left(\ln{4L_fR^2} + \ln{\varepsilon^{-1}} \right) = O\left(\sqrt{\frac{L_f}{\mu}} \ln{\frac{L_f R^2}{\varepsilon}} \right),
        \end{gathered}
    \end{equation*}
we guarantee that 
\begin{equation*}
        f(x_N) - f(x^*) \leqslant \varepsilon. 
\end{equation*}

\end{remark}

\begin{remark}
    \label{remark 5.5}
    Let us consider the setting of Remark~\ref{remark4.2}, where we made a reduction of the convex case $\mu=0$ to the strongly convex case by introducing a quadratic regularization with regularization parameter $\mu$. Recall that this led to the bound
     \begin{multline*}
        f(x_N) - f(x^*) \leqslant L R^2 \exp\left({-\frac{1}{2}\sqrt{\frac{\mu}{2(L + 2)}}N} \right) + \left( 1 + \sqrt{\frac{2L + 4}{\mu}} \right) \left( \frac{1}{L} + \frac{1}{\mu} \right) \delta^2 + \frac{\mu}{2}R^2,
    \end{multline*}
    where $R$ is such that $\|\tilde{x}_0-x^*\|_2 \leq R$.
    We choose the regularization parameter $\mu$, the error level $\delta$, and the number of iterations $N$ such that each of the three terms in this bound are smaller than $\frac{\varepsilon}{3}$. Then, choosing
     \begin{equation*}
        \mu = \frac{2}{3} \frac{\varepsilon}{R^2},
    \end{equation*}
    \begin{equation*}
        \begin{gathered}
            \delta \leqslant \left(\frac{2}{243} \right)^{\frac{1}{4}} \frac{1}{\sqrt{1 + \sqrt{2L + 4}}} R^{-\frac{3}{2}} \varepsilon^{\frac{5}{4}} = O \left(L^{-\frac{1}{4}} R^{-\frac{3}{2}} \varepsilon^{\frac{5}{4}} \right),
        \end{gathered}
    \end{equation*}
    \begin{equation*}
        \begin{gathered}
            N \geqslant \sqrt{12L + 24}R\ln{2LR^2} + 2\sqrt{2L + 4}\frac{1}{\sqrt{\varepsilon}}\ln{\frac{1}{\varepsilon}} = O\left(\sqrt{\frac{LR^2}{\varepsilon}} \ln{\frac{LR^2}{\varepsilon}} \right),
        \end{gathered}
    \end{equation*}
    we guarantee that 
\begin{equation*}
        f(x_N) - f(x^*) \leqslant \varepsilon. 
\end{equation*}

\end{remark}

\begin{remark}
    \label{remark5.6}
    Let us apply Theorem~\ref{stopping_rule} for solving linear inverse problems. Let $A \in \mathbb{R}^{n \times n} $ be such that $\det(A) \not= 0$ and consider the following linear system for finding $x \mathbb{R}^{n}$: $Ax=b$.
    Solving this problem is equivalent to solving the convex optimization problem:
    \begin{equation*}
        \begin{gathered}
           \min_{x \in \mathbb{R}^n} \left\{ f(x) = \frac{1}{2} \|Ax - b \|_2^2 \right\}.
        \end{gathered}
    \end{equation*}
    If we solve the latter problem with accuracy $\varepsilon=\frac{\varepsilon_0^2}{2}$, then we guarantee that $\|Ax - b\|_2 \leqslant \varepsilon_0$. 
    
    %
   Let us assume that the solution $x^*$ satisfies $\|x^* \|_2 \leqslant R_*$ and that Algorithm \ref{alg3} starts from the point $0$. Then, we can take $R=R_*$.
    %
    %
    According to Theorem~\ref{stopping_rule}, given a target accuracy $\zeta>0$, we have that Algorithm \ref{alg3} stops after $N_{\pd{\text{stop}}}$ iterations such that  $N_{\pd{\text{stop}}} \leq  \sqrt{\frac{2LR_*^2}{\zeta}}  + 1$. Moreover, we have that
        \begin{equation*}
        \begin{gathered}
            f(x_{N_{\text{stop}}}) - f(x^*) \leqslant 
            \frac{\delta_2}{A_N} \sum_{j=0}^{N_{\text{stop}}} A_j + 3\delta_1 R_* + \zeta   
            \leqslant  N_{\pd{\text{stop}}} \delta_2 + 3\delta_1 R_* + \zeta 
        \end{gathered}
    \end{equation*}  
    since $\lbrace A_j \rbrace$ is an increasing sequence.
    
    Choosing $\zeta  \leqslant \frac{\varepsilon}{3}$ and $\delta \leqslant \min \left\{\left(\frac{L^{\frac{1}{4}}}{6 \sqrt{3R_*}} \right) \varepsilon^{\frac{3}{4}},   \frac{\varepsilon}{9R_*}\right\}$, we guarantee that $f(x_{N_{\text{stop}}}) - f(x^*) \leqslant \varepsilon$ and, hence, $\|Ax - b\|_2 \leqslant \varepsilon_0$. Moreover, the number of iterations to guarantee such a solution is bounded as
    \begin{equation*}
        N_{\varepsilon_0} = \frac{\sqrt{6LR_*^2}}{\varepsilon_0} + 1.
    \end{equation*}
    
\end{remark}

\begin{remark}
\label{garnot compare}
In the setting of relative noise in the gradient, Theorem \ref{th_alph} says that whenever $\alpha \leqslant O \left( \frac{\mu}{L} \right)$, STM converges linearly in the same way as accelerated gradient method in the exact setting, i.e., with the rate $O\left(\left(1-\sqrt{\frac{\mu}{L}}\right)^k \right)$ which is faster than the convergence rate $O\left(\left(1-\frac{\mu}{L}\right)^k \right)$ of gradient descent. Here $k$ is the iteration counter. The paper \cite{gannot2021frequency} considers, in particular, accelerated method, called the  Triple Momentum Method, in the presence of relative noise in the gradient. They show that when $\alpha < \frac{\sqrt{\chi} + 1}{4 \chi - 3 \sqrt{\chi} + 1} = O \left( \sqrt{\frac{\mu}{L}} \right)$, where $\chi = \frac{L}{\mu}$, the   Triple Momentum Method converges with a linear rate as well. At the same time, their convergence rate depends on the noise level $\alpha$, is no better than the accelerated rate, and is equal to it only in the case $\alpha=0$.
Figure \ref{figure_tmm_rate} illustrates the situation for two different values of the condition number $\chi$. The black dashed line shows the threshold, below which STM with relative inexactness in the gradient has linear convergence rate similar to exact STM, and the latter rate is denoted by the orange line. Green line shows the convergence rate of the gradient method. Finally, the blue line shows the dependence of the convergence rate in \cite{gannot2021frequency} on the inexactness level $\alpha$. As we see, it can be even worse than that of the gradient method for large values of $\alpha$.

\begin{figure}[H] 
	\center{\includegraphics[scale=0.35]{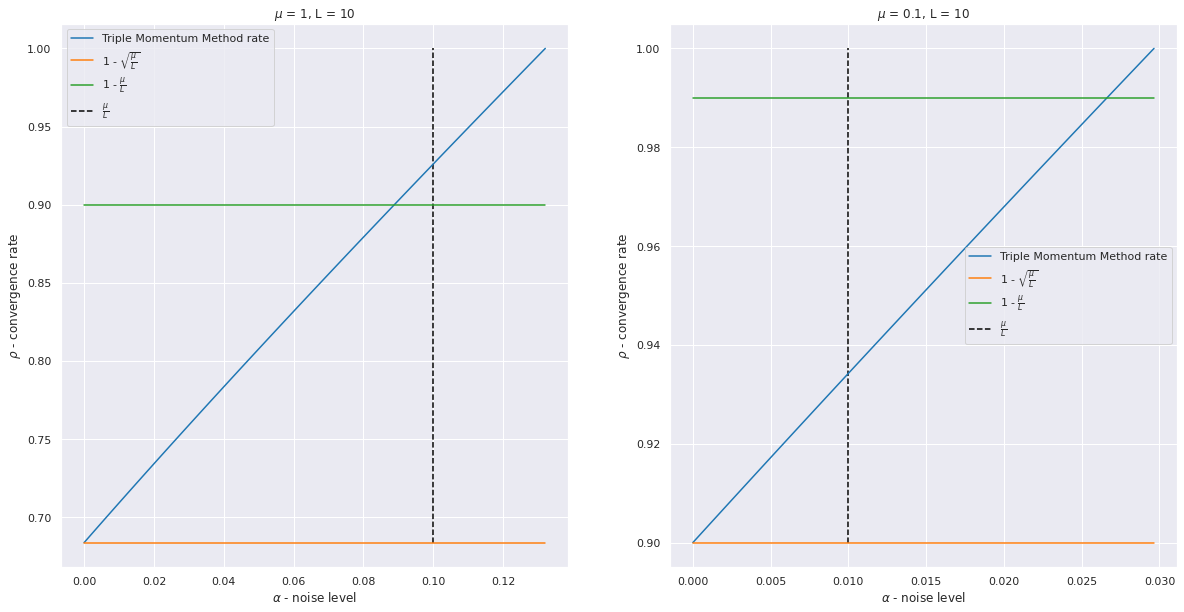}}
	\caption{Comparison of the convergence rate of Triple Momentum Method and STM}
	\label{figure_tmm_rate}
\end{figure}

As our experiments show, STM is more robust in the relative noise setting, that is, numerically estimating the dependence of the largest possible  $\alpha = \alpha^*$ for given problem parameters $\mu, L$, we get a larger upper bound. More detailed information \pd{can be found} in Section~\ref{section exp}. This leads us to the hypothesis that the condition $\alpha \leqslant O \left( \frac{\mu}{L} \right)$ for inexact STM may be weakened.
\end{remark}

\section{Numerical experiments}
\label{section exp}
In this section, we provide a series of numerical experiments to illustrate the practical performance of the considered algorithms under absolute and relative noise. The noise was generated as independent random uniform and unbiased.


We start with the experiments in the setting of $\mu=0$ using the following objective function described in \cite[p. 69]{nesterov2018lectures} and known as the worst-case function for first-order methods:
\begin{equation*}
    \begin{gathered}
        f(x) = \frac{L}{8} \left( x_1^2 + \displaystyle\sum_{j = 0}^{k - 1} \left(x_j - x_{j + 1}\right)^2 + x_k^2 \right) - \frac{L}{4} x_1, \\
        x^* = \left(1 - \frac{1}{k + 1}, \: \dots \:,  1 - \frac{k}{k + 1}, \: 0, \: \dots \:, 0\right)^T, \\
        1 \leqslant k \leqslant \dim{x}.
    \end{gathered}
\end{equation*}
The next two plots show the convergence of STM at the first 50 000 and 10 000 iterations, respectively, in the absolute noise setting with different values of $\delta$.
%
\begin{figure}[H]
	\center{\includegraphics[width=1\linewidth]{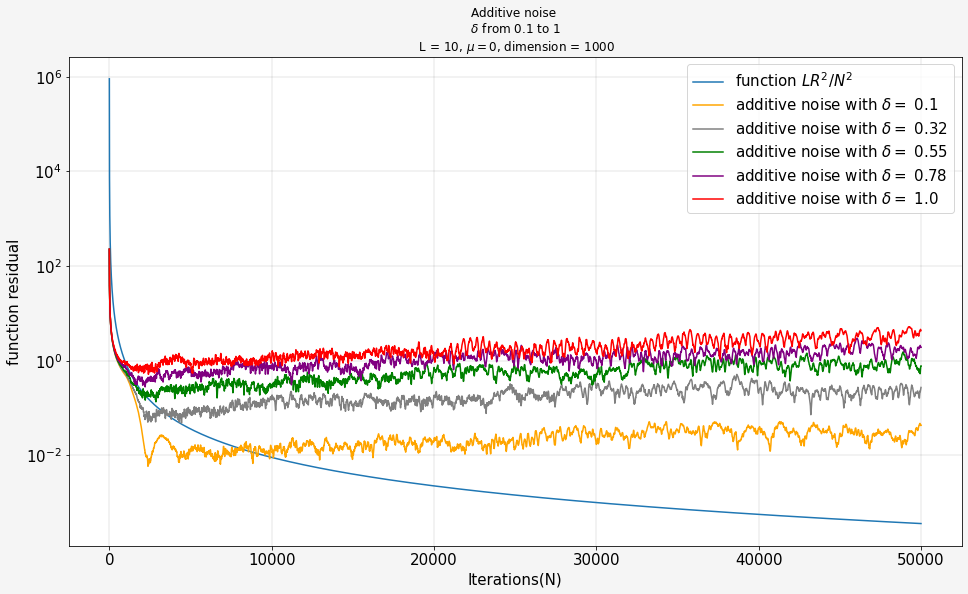}}
    \caption{First test -- the performance of STM for $\mu = 0$ for the first 50 000 iterations.}
\end{figure}
\begin{figure}[H]
	\center{\includegraphics[width=1\linewidth]{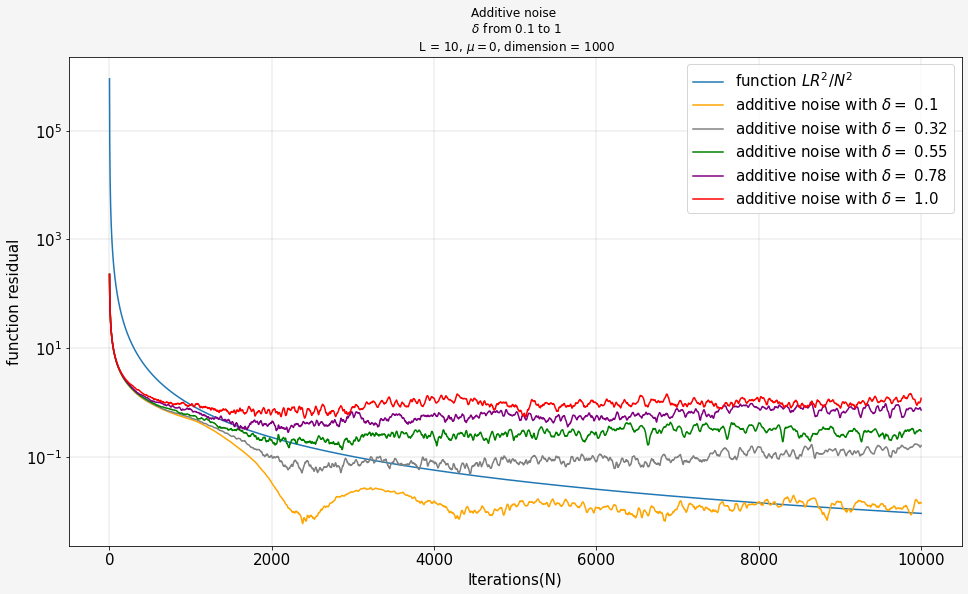}}
	\caption{First test -- the performance of STM for $\mu = 0$ for the first 10 000 iterations.}
\end{figure}
We can observe that, as predicted by Theorem \ref{Th:STM_main}, we see that the increasing third term in the convergence rate \eqref{eq:STM_convergence_conv} at some point starts to overweight the first decreasing term.  

We further compare the convergence in two different settings of the noise: absolute and relative.
\begin{figure}[H]
	\begin{center}
		\includegraphics[width=\textwidth]{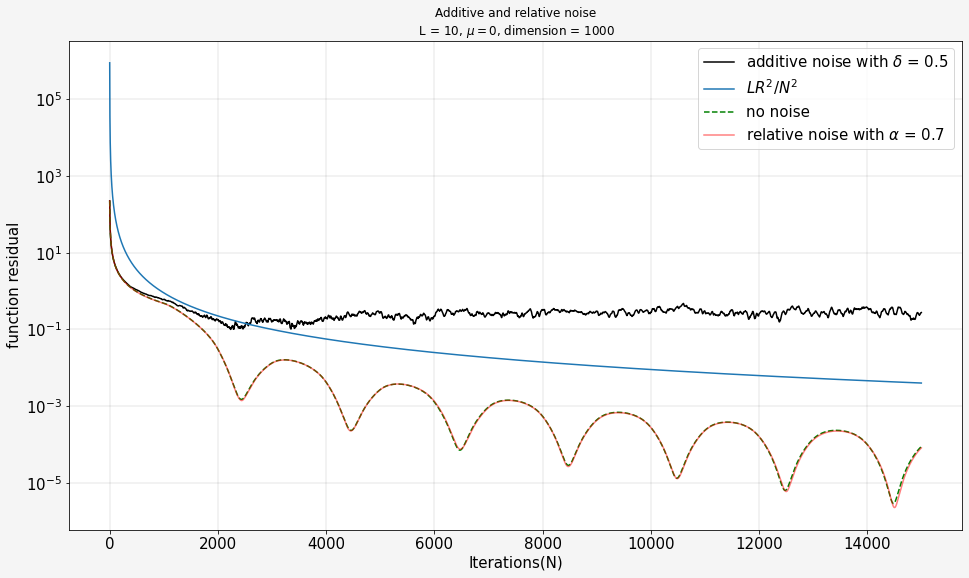}
	\caption{Second test -- the performance of STM for $\mu = 0$ with relative and additive types of noise.}
	\end{center}
\end{figure}
\pd{
As was expected from the theory, for sufficiently small $\alpha$, the convergence of inexact method is very close to the convergence of the exact method. 
Since in this experiment the noise is stochastic, this effect can be possibly explained using the theoretical results obtained in \cite{vaswani2019fast}: under the strong growth condition (SGC)
\begin{equation*}
    \mathbb{E}_{\xi} \| \tilde{\nabla} f(x, \xi) \|_2^2 \leqslant \rho \|\nabla f(x) \|_2^2,
\end{equation*}
$L_f$-smoothness and convexity, SGD with Nesterov's acceleration has the following convergence rate:
\begin{equation*}
    \mathbb{E} f(x^k) - f(x^*) \leqslant \frac{2\rho^2 L_f}{k^2} \|x_0 - x^* \|_2^2,
\end{equation*}
i.e., similar to the deterministic method despite that the gradients are stochastic.
SGC can be translated into the relative noise condition~\eqref{relative_inexact}, making them related.
Although a different method is used in our paper, the obtained results make it reasonable to expect a similar convergence in the concept of relative noise as in the absence of any noise.
}

The next plot illustrates the convergence of STM in the setting of $\mu=0$ and relative noise in the gradient for different values of the parameter $\alpha$.
\begin{figure}[H]
	\begin{center}
		\includegraphics[width=0.8 \textwidth]{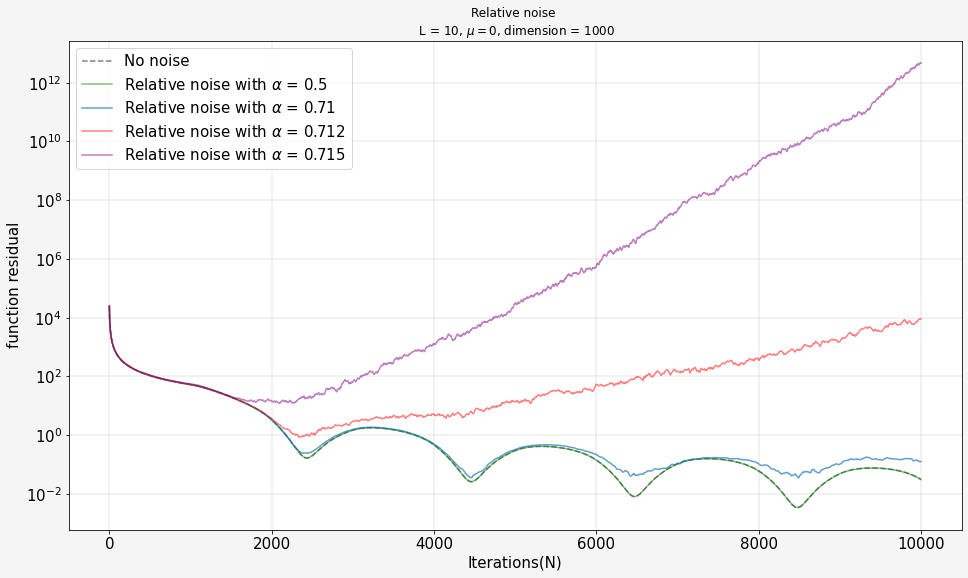}
	\caption{Third test -- the performance of STM with relative noise and $\mu = 0$ for different values of $\alpha$.}
	\end{center}
\end{figure}
As we see, for $\alpha \leqslant 0.71$, the convergence of the method does not deteriorate and the value $\alpha^* \approx 0.71$ can be seen as a threshold above which the method diverges.

We next explore the strongly convex setting with $\mu>0$ using the worst-case function \cite[p.78]{nesterov2018lectures}:
\begin{equation*}
    \begin{gathered}
        f(x) = \frac{\mu \left(\chi - 1\right)}{8} \left(x_1^2 + \displaystyle\sum_{j = 1}^{n - 1} {\left(x_j - x_{j + 1}\right)^2} - 2x_1\right) + \frac{\mu}{2} \|x \|_2^2, \\
        \chi = \frac{L}{\mu}. 
    \end{gathered}
\end{equation*}

We first consider the performance of STM with absolute noise for different values of $\delta$. Dashed lines represent the corresponding theoretical bound.
\begin{figure}[H]
	\center{\includegraphics[width=0.8\linewidth]{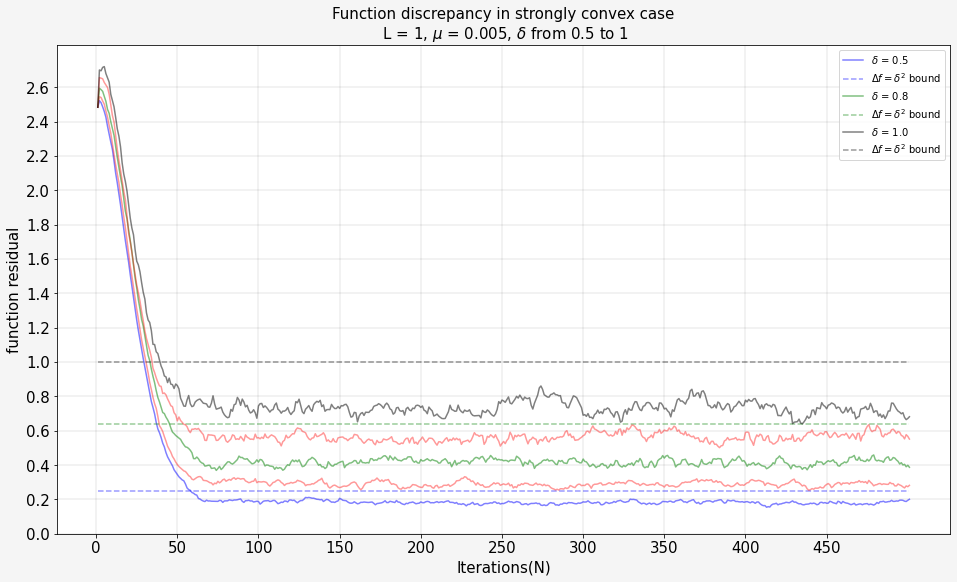}}
    \caption{Fourth test -- the performance of STM for $\mu > 0$ and absolute noise $\delta \in \lbrace 0.5, 0.6, 0.7, 0.8, 0.9, 1.0 \rbrace.$}
\end{figure}
The next plot confirms Theorem~\ref{th4conv} and Remark~\ref{remark5.4}.
\begin{figure}[H]
	\center{\includegraphics[width=1\linewidth]{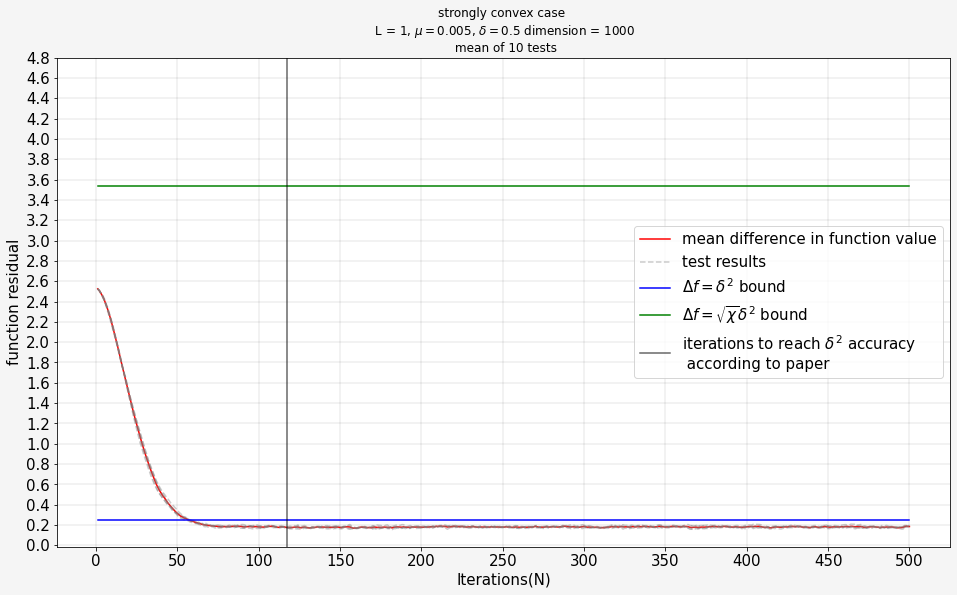}}
	\caption{Fifth test -- mean of 30 tests, level of approximation and required number of steps.}
\end{figure}
Next, similarly to the degenerate case $\mu=0$, we consider the behavior of the method for different parameters $\alpha$ when a relative noise is present in the gradient.
\begin{figure}[H]
	\center{\includegraphics[width=0.9\linewidth]{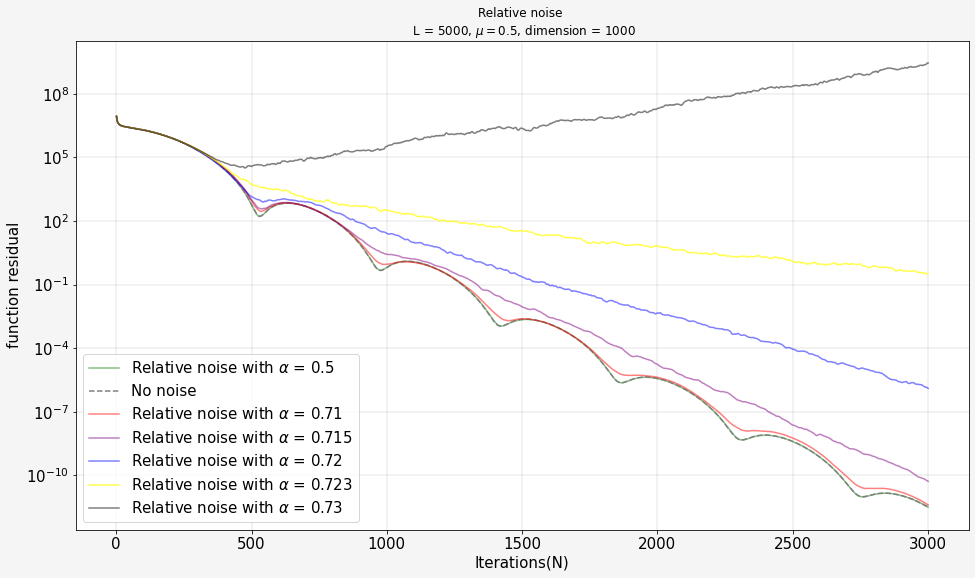}}
	\caption{Sixth test -- the performance of STM with relative noise and $\mu > 0$ for different values of $\alpha$.}
\end{figure}
Note that, in the strongly convex case, we observe a similar effect as in the degenerate case:  the algorithm converges for $\alpha$-values smaller than a certain threshold value $\alpha^*$.

\begin{figure}[H] 
	\center{\includegraphics[width=1\linewidth]{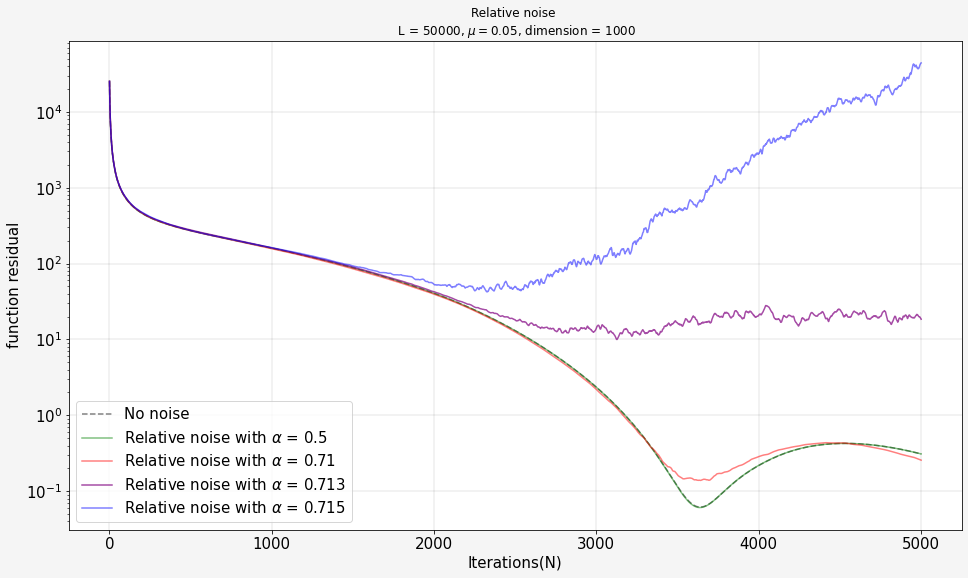}}
	\caption{Seventh test -- the performance of STM with relative noise and $\mu > 0$ for different values of $\alpha$.}
\end{figure}

\begin{figure}[H] 
	\center{\includegraphics[width=0.9\linewidth]{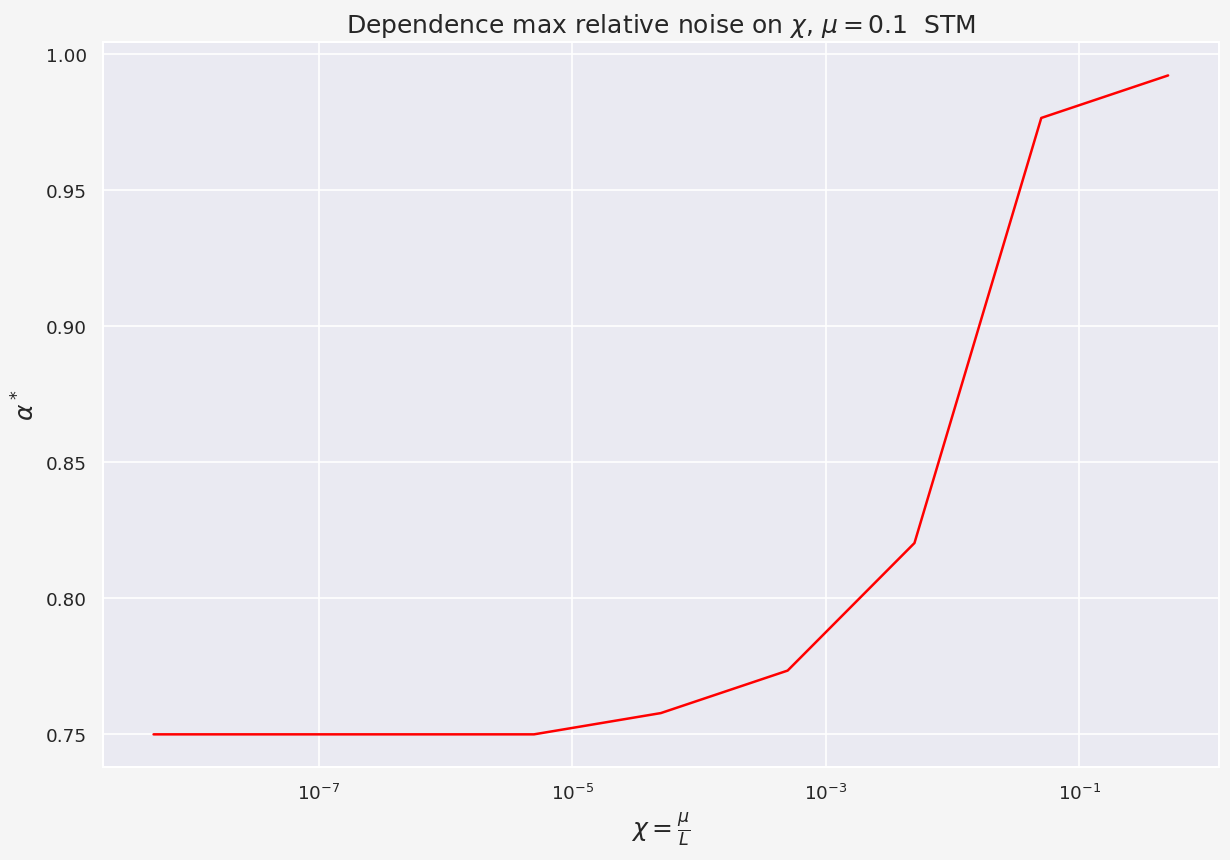}}
	\caption{Eigth test -- threshold $\alpha^*$ for different $L$ and $\mu = 0.1$, for STM}
	\label{figure_rel_stm}
\end{figure}

\begin{figure}[H] 
	\center{\includegraphics[width=0.9\linewidth]{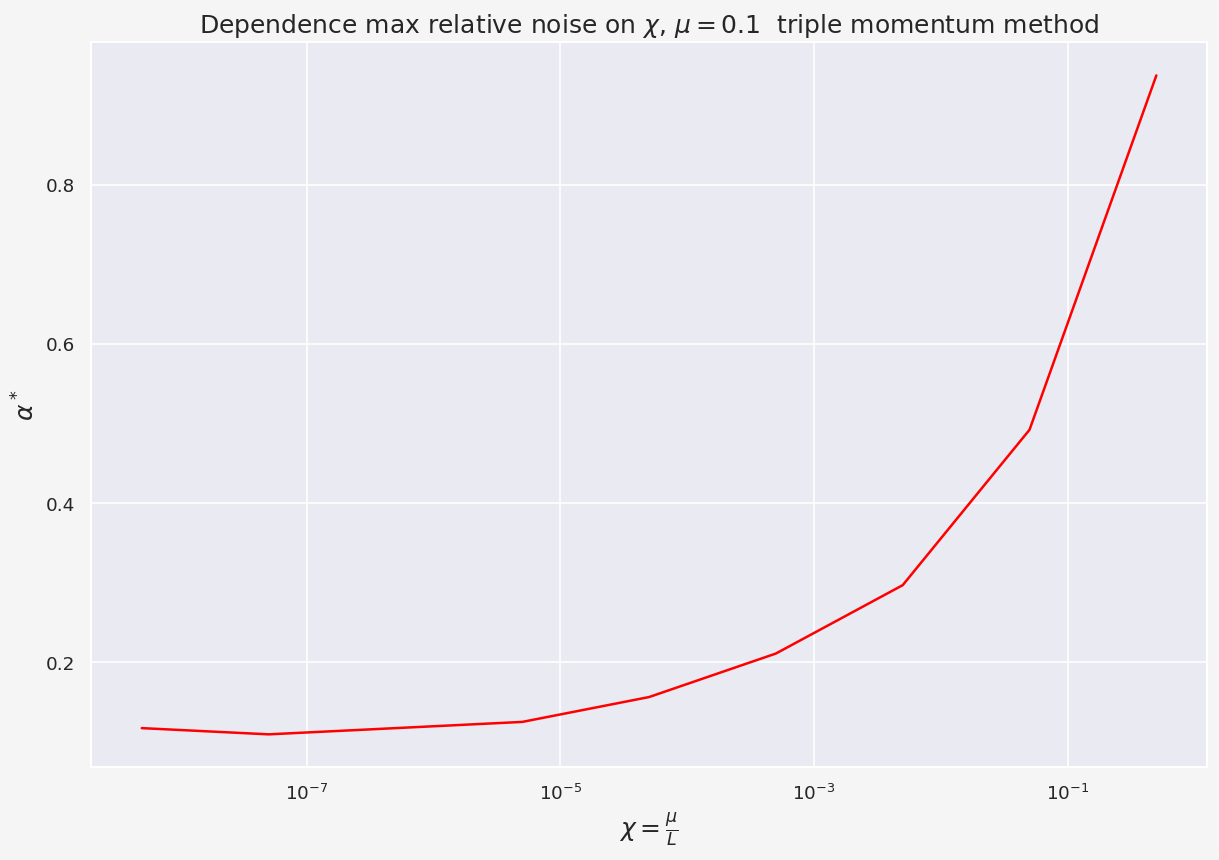}}
	\caption{Nineth test -- threshold $\alpha^*$ for different $L$ and $\mu = 0.1$, for the Triple Momentum Method}
	\label{figure_rel_trip}
\end{figure}

Finally, we compare STM and triple momentum method.
Figures~$\ref{figure_rel_stm}$ and $\ref{figure_rel_trip}$ show, that for the same parameters of the problem, STM is capable of converging at a much higher noise level than triple momentum algorithm.

%
\section*{Acknowledgments}
The authors are grateful to Eduard Gorbunov for useful discussions.

\bibliographystyle{tfs}
\bibliography{literature,all_refs3}

\begin{thebibliography}{10}
\providecommand{\MR}{\relax\unskip\space MR }
\providecommand{\url}[1]{\normalfont{#1}}
\providecommand{\urlprefix}{Available at }

\bibitem{ajalloeian2020analysis}
A. Ajalloeian and S.U. Stich, \emph{Analysis of sgd with biased gradient
  estimators}, deepai  (2020).

\bibitem{akhavan2020exploiting}
A. Akhavan, M. Pontil, and A. Tsybakov, \emph{Exploiting higher order
  smoothness in derivative-free optimization and continuous bandits}, Advances
  in Neural Information Processing Systems 33 (2020), pp. 9017--9027.

\bibitem{bach16highly-smooth}
F. Bach and V. Perchet, \emph{Highly-Smooth Zero-th Order Online Optimization},
  in \emph{29th Annual Conference on Learning Theory}, V. Feldman, A. Rakhlin,
  and O. Shamir, eds., Proceedings of Machine Learning Research Vol.~49, 23--26
  Jun, Columbia University, New York, New York, USA. PMLR, 2016, pp. 257--283.
  \urlprefix\url{http://proceedings.mlr.press/v49/bach16.html}.

\bibitem{beck2017first}
A. Beck, \emph{First-order methods in optimization}, SIAM, 2017.

\bibitem{belloni2015escaping}
A. Belloni, T. Liang, H. Narayanan, and A. Rakhlin, \emph{Escaping the Local
  Minima via Simulated Annealing: Optimization of Approximately Convex
  Functions}, in \emph{Proceedings of The 28th Conference on Learning Theory},
  P. Grünwald, E. Hazan, and S. Kale, eds., Proceedings of Machine Learning
  Research Vol.~40, 03--06 Jul, Paris, France. PMLR, 2015, pp. 240--265.
  \urlprefix\url{http://proceedings.mlr.press/v40/Belloni15.html}.

\bibitem{ben-tal2015lectures}
A. Ben-Tal and A. Nemirovski, \emph{Lectures on Modern Convex Optimization
  (Lecture Notes)}, Personal web-page of A. Nemirovski, 2015.

\bibitem{berahas2021theoretical}
A.S. Berahas, L. Cao, K. Choromanski, and K. Scheinberg, \emph{A theoretical
  and empirical comparison of gradient approximations in derivative-free
  optimization}, Foundations of Computational Mathematics  (2021), pp. 1--54.

\bibitem{beznosikov2020gradient}
A. Beznosikov, A. Sadiev, and A. Gasnikov, \emph{Gradient-Free Methods with
  Inexact Oracle for Convex-Concave Stochastic Saddle-Point Problem}, in
  \emph{International Conference on Mathematical Optimization Theory and
  Operations Research}. Springer, 2020, pp. 105--119.

\bibitem{bubeck2015convex}
S. Bubeck, \emph{et~al.}, \emph{Convex optimization: Algorithms and
  complexity}, Foundations and Trends{\textregistered} in Machine Learning 8
  (2015), pp. 231--357.

\bibitem{cohen2018acceleration}
M. Cohen, J. Diakonikolas, and L. Orecchia, \emph{On acceleration with
  noise-corrupted gradients}, in \emph{International Conference on Machine
  Learning}. PMLR, 2018, pp. 1019--1028.

\bibitem{conn2009introduction}
A. Conn, K. Scheinberg, and L. Vicente, \emph{Introduction to Derivative-Free
  Optimization}, Society for Industrial and Applied Mathematics, 2009,
  \urlprefix\url{http://epubs.siam.org/doi/abs/10.1137/1.9780898718768}.

\bibitem{d2008smooth}
A. d'Aspremont, \emph{Smooth optimization with approximate gradient}, SIAM
  Journal on Optimization 19 (2008), pp. 1171--1183.

\bibitem{devolder2011stochastic}
O. Devolder, \emph{Stochastic first order methods in smooth convex
  optimization}, CORE Discussion Paper 2011/70  (2011).

\bibitem{devolder2013exactness}
O. Devolder, \emph{Exactness, inexactness and stochasticity in first-order
  methods for large-scale convex optimization}, Ph.D. diss., ICTEAM and CORE,
  Universit{\'e} Catholique de Louvain,  2013.

\bibitem{devolder2014first}
O. Devolder, F. Glineur, and Y. Nesterov, \emph{First-order methods of smooth
  convex optimization with inexact oracle}, Mathematical Programming 146
  (2014), pp. 37--75.
  \urlprefix\url{http://dx.doi.org/10.1007/s10107-013-0677-5}.

\bibitem{devolder2013first}
O. Devolder, F. Glineur, Y. Nesterov, \emph{et~al.}, \emph{First-order methods
  with inexact oracle: the strongly convex case}, CORE Discussion Papers
  2013016 (2013), p.~47.

\bibitem{drusvyatskiy2022stochastic}
D. Drusvyatskiy and L. Xiao, \emph{Stochastic optimization with
  decision-dependent distributions}, Mathematics of Operations Research
  (2022).

\bibitem{dvinskikh2019decentralized}
D. Dvinskikh and A. Gasnikov, \emph{Decentralized and parallelized primal and
  dual accelerated methods for stochastic convex programming problems}, Journal
  of Inverse and Ill-posed Problems  (2021).

\bibitem{dvinskikh2020accelerated}
D.M. Dvinskikh, A.I. Turin, A.V. Gasnikov, and S.S. Omelchenko,
  \emph{Accelerated and non accelerated stochastic gradient descent in model
  generality}, Matematicheskie Zametki 108 (2020), pp. 515--528.

\bibitem{dvurechensky2020numerical}
P. Dvurechensky, \emph{Numerical methods in large-scale optimization: inexact
  oracle and primal-dual analysis}, HSE. Habilitation  (2020).

\bibitem{dvurechensky2016stochastic}
P. Dvurechensky and A. Gasnikov, \emph{Stochastic intermediate gradient method
  for convex problems with stochastic inexact oracle}, Journal of Optimization
  Theory and Applications 171 (2016), pp. 121--145.
  \urlprefix\url{http://dx.doi.org/10.1007/s10957-016-0999-6}.

\bibitem{dvurechensky2018computational}
P. Dvurechensky, A. Gasnikov, and A. Kroshnin, \emph{Computational Optimal
  Transport: Complexity by Accelerated Gradient Descent Is Better Than by
  {S}inkhorn’s Algorithm}, in \emph{Proceedings of the 35th International
  Conference on Machine Learning}, J. Dy and A. Krause, eds., Proceedings of
  Machine Learning Research Vol.~80. 2018, pp. 1367--1376. arXiv:1802.04367.

\bibitem{dvurechensky2021first}
P. Dvurechensky, S. Shtern, and M. Staudigl, \emph{First-order methods for
  convex optimization}, EURO Journal on Computational Optimization 9 (2021), p.
  100015.

\bibitem{d2021acceleration}
A. d’Aspremont, D. Scieur, A. Taylor, \emph{et~al.}, \emph{Acceleration
  methods}, Foundations and Trends{\textregistered} in Optimization 5 (2021),
  pp. 1--245.

\bibitem{evtushenko2013optimization}
Y.G. Evtushenko, \emph{Optimization and fast automatic differentiation},
  Computing Center of RAS, Moscow  (2013).

\bibitem{gannot2021frequency}
O. Gannot, \emph{A frequency-domain analysis of inexact gradient methods},
  Mathematical Programming 194 (2022), pp. 975--1016.
  \urlprefix\url{https://doi.org/10.1007/s10107-021-01665-8}.

\bibitem{gasnikov2016efficient}
A.V. Gasnikov, E.V. Gasnikova, Y.E. Nesterov, and A.V. Chernov, \emph{Efficient
  numerical methods for entropy-linear programming problems}, Computational
  Mathematics and Mathematical Physics 56 (2016), pp. 514--524.
  \urlprefix\url{http://dx.doi.org/10.1134/S0965542516040084}.

\bibitem{gasnikov2017modern}
A. Gasnikov, \emph{Universal gradient descent}, arXiv preprint arXiv:1711.00394
   (2017).

\bibitem{gasnikov2017convex}
A. Gasnikov, S. Kabanikhin, A. Mohammed, and M. Shishlenin, \emph{Convex
  optimization in hilbert space with applications to inverse problems}, arXiv
  preprint arXiv:1703.00267  (2017).

\bibitem{gasnikov2018universal}
A.V. Gasnikov and Y.E. Nesterov, \emph{Universal method for stochastic
  composite optimization problems}, Computational Mathematics and Mathematical
  Physics 58 (2018), pp. 48--64.

\bibitem{goodfellow2016deep}
I. Goodfellow, Y. Bengio, A. Courville, and Y. Bengio, \emph{Deep learning},
  Vol.~1, MIT press Cambridge, 2016.

\bibitem{gorbunov2019optimal}
E. Gorbunov, D. Dvinskikh, and A. Gasnikov, \emph{Optimal decentralized
  distributed algorithms for stochastic convex optimization}, arXiv preprint
  arXiv:1911.07363  (2019).

\bibitem{gorbunov2018accelerated}
E. Gorbunov, P. Dvurechensky, and A. Gasnikov, \emph{An accelerated method for
  derivative-free smooth stochastic convex optimization}, arXiv preprint
  arXiv:1802.09022  (2018).

\bibitem{kabanikhin2011inverse}
S.I. Kabanikhin, \emph{Inverse and ill-posed problems: theory and
  applications}, Vol.~55, Walter De Gruyter, 2011.

\bibitem{kamzolov2020universal}
D. Kamzolov, P. Dvurechensky, and A.V. Gasnikov, \emph{Universal intermediate
  gradient method for convex problems with inexact oracle}, Optimization
  Methods and Software  (2020), pp. 1--28.

\bibitem{kotsalis2022simple}
G. Kotsalis, G. Lan, and T. Li, \emph{Simple and optimal methods for stochastic
  variational inequalities, i: operator extrapolation}, SIAM Journal on
  Optimization 32 (2022), pp. 2041--2073.

\bibitem{lan2020first}
G. Lan, \emph{First-order and Stochastic Optimization Methods for Machine
  Learning}, Springer, 2020.

\bibitem{nemirovski1986regularizing}
A.S. Nemirovski, \emph{Regularizing properties of the conjugate gradient method
  for ill-posed problems}, Zhurnal Vychislitel'noi Matematiki i Matematicheskoi
  Fiziki 26 (1986), pp. 332--347.

\bibitem{nemirovsky1983problem}
A. Nemirovsky and D. Yudin, \emph{Problem Complexity and Method Efficiency in
  Optimization}, J. Wiley \& Sons, New York, 1983.

\bibitem{nesterov2018lectures}
Y. Nesterov, \emph{Lectures on convex optimization}, Vol. 137, Springer, 2018.

\bibitem{nesterov2017random}
Y. Nesterov and V. Spokoiny, \emph{Random gradient-free minimization of convex
  functions}, Found. Comput. Math. 17 (2017), pp. 527--566.
  \urlprefix\url{https://doi.org/10.1007/s10208-015-9296-2}, First appeared in
  2011 as CORE discussion paper 2011/16.

\bibitem{novitskii2021improved}
V. Novitskii and A. Gasnikov, \emph{Improved exploiting higher order smoothness
  in derivative-free optimization and continuous bandit}, arXiv preprint
  arXiv:2101.03821  (2021).

\bibitem{pedregosa2020average}
F. Pedregosa and D. Scieur, \emph{Average-case acceleration through spectral
  density estimation}, arXiv preprint arXiv:2002.04756  (2020).

\bibitem{poljak1981iterative}
B. Poljak, \emph{Iterative algorithms for singular minimization problems}, in
  \emph{Nonlinear Programming 4}, Elsevier,  1981, pp. 147--166.

\bibitem{polyak1987introduction}
B. Polyak, \emph{Introduction to Optimization}, New York, Optimization
  Software, 1987.

\bibitem{polyak1990optimal}
B.T. Polyak and A.B. Tsybakov, \emph{Optimal order of accuracy of search
  algorithms in stochastic optimization}, Problemy Peredachi Informatsii 26
  (1990), pp. 45--53.

\bibitem{risteski2016algorithms}
A. Risteski and Y. Li, \emph{Algorithms and matching lower bounds for
  approximately-convex optimization}, Advances in Neural Information Processing
  Systems 29 (2016), pp. 4745--4753.

\bibitem{rockafellar1970convex}
R.T. Rockafellar, \emph{Convex analysis}, Vol.~36, Princeton university press,
  1970.

\bibitem{scieur2020universal}
D. Scieur and F. Pedregosa, \emph{Universal Asymptotic Optimality of Polyak
  Momentum}, in \emph{International Conference on Machine Learning}. PMLR,
  2020, pp. 8565--8572.

\bibitem{stonyakin2021inexact}
F. Stonyakin, A. Tyurin, A. Gasnikov, P. Dvurechensky, A. Agafonov, D.
  Dvinskikh, M. Alkousa, D. Pasechnyuk, S. Artamonov, and V. Piskunova,
  \emph{Inexact model: A framework for optimization and variational
  inequalities}, Optimization Methods and Software  (2021).
  \urlprefix\url{https://doi.org/10.1080/10556788.2021.1924714}.

\bibitem{stonyakin2020adaptive}
F. Stonyakin, \emph{Adaptive methods for variational inequalities, minimization
  problems and functional with generalized growth condition}, MIPT.
  Habilitation  (2020).

\bibitem{taylor2017smooth}
A.B. Taylor, J.M. Hendrickx, and F. Glineur, \emph{Smooth strongly convex
  interpolation and exact worst-case performance of first-order methods},
  Mathematical Programming 161 (2017), pp. 307--345.

\bibitem{tyurin2020adaptive}
A. Tyurin, \emph{Development of a method for solving structural optimization
  problems}, HSE. PhD Thesis  (2020).

\bibitem{vasilyev2002optimization}
F. Vasilyev, \emph{Optimization Methods}, Moscow, Russia: FP, 2002.

\bibitem{vaswani2019fast}
S. Vaswani, F. Bach, and M. Schmidt, \emph{Fast and faster convergence of sgd
  for over-parameterized models and an accelerated perceptron}, in \emph{The
  22nd International Conference on Artificial Intelligence and Statistics}.
  PMLR, 2019, pp. 1195--1204.

\end{thebibliography}
\appendix

\end{document}